\documentclass[a4paper,11pt]{amsart}

\usepackage {amssymb}
\usepackage{amscd}

\theoremstyle{plain}
\newtheorem{thm}{Theorem}[section]
\newtheorem{lemma}[thm]{Lemma}
\newtheorem{cor}[thm]{Corollary}
\newtheorem{prop}[thm]{Proposition}

\theoremstyle{definition}

\newtheorem{rmk}[thm]{Remark}
\newtheorem{example}[thm]{Example}
\newtheorem{examples}[thm]{Examples}

\def\Gr{\operatorname{Gr}}

\def\ad{\operatorname{ad}}
\def\deg{\operatorname{deg}}
\def\det{\operatorname{det}}

\def\hom{\operatorname{Hom}}
\def\Hom{\operatorname{Hom}}
\def\ker{\operatorname{Ker}}
\def\im{\operatorname{Im}}

\def\Ker{\operatorname{Ker}}
\def\Im{\operatorname{Im}}

\def\Id{\operatorname{Id}}
\def\deg{\operatorname{deg}}

\def\top{\operatorname{top}}

\newcommand{\field}{\mathbb}

\newcommand{\C}{{\field C}}

\newcommand{\hra}{\hookrightarrow}

\newcommand{\frb}{\mathfrak{b}}

\newcommand{\frg}{\mathfrak{g}}
\newcommand{\frh}{\mathfrak{h}}

\newcommand{\frk}{\mathfrak{k}}
\newcommand{\frl}{\mathfrak{l}}

\newcommand{\frn}{\mathfrak{n}}
\newcommand{\fro}{\mathfrak{o}}
\newcommand{\frp}{\mathfrak{p}}
\newcommand{\frqq}{\mathfrak{q}}
\newcommand{\frr}{\mathfrak{r}}
\newcommand{\frs}{\mathfrak{s}}
\newcommand{\frt}{\mathfrak{t}}
\newcommand{\fru}{\mathfrak{u}}

\newcommand{\bbC}{\mathbb{C}}

\newcommand{\bbN}{\mathbb{N}}

\newcommand{\bbR}{\mathbb{R}}

\newcommand{\bbZ}{\mathbb{Z}}

\newcommand{\caD}{\mathcal{D}}

\newcommand{\del}{\partial}
\newcommand{\no}{\noindent}

\begin{document}

\title{Dirac operators and Lie algebra cohomology}

\author{Jing-Song Huang}
\address{Department of Mathematics, Hong Kong University of Science and Technology, Clear Water Bay, Kowloon, Hong Kong SAR, China}
\email{mahuang\char'100ust.hk}

\author{Pavle Pand\v zi\'c} 
\address{Department of Mathematics, University of Zagreb, Bijeni\v cka 30,
10000 Zagreb, Croatia}
\email{pandzic\char'100math.hr}

\author{David Renard}
\address{Centre de math\'ematiques Laurent Schwartz, Ecole Polytechnique, 91128 Palaiseau Cedex, France}
\email{renard@math.polytechnique.fr}

\abstract{Dirac cohomology is a new tool to study unitary and admissible representations
of semisimple Lie groups. It was introduced by Vogan and further studied by Kostant and
ourselves \cite{V2}, \cite{HP1}, \cite{Kdircoh}.
The aim of this paper is to study the Dirac cohomology for the 
Kostant cubic Dirac operator and its relation to Lie algebra cohomology.
We show that the Dirac cohomology coincides with the corresponding nilpotent 
Lie algebra cohomology in many cases,
while in general it has better algebraic behavior and it is more accessible for calculation.}\endabstract


\thanks{The research of the first named author was partially supported by
RGC-CERG grants of Hong Kong SAR and National Nature Science Foundation of China.
The research of the second named author was partially supported by a grant
from the Ministry of Science 
and Technology of Republic of Croatia. The second and third named
authors were also supported by the joint French-Croatian program
COGITO. 
Parts of this work were done during authors' visits to CNRS,  University of Paris VII,
Mathematisches Forschungsinstitut Oberwolfach, and Institute of Mathematical Sciences and Department of Mathematics at
the National  University of Singapore.  The authors thank these institutions for their generous support and hospitality.}


\maketitle

\section{Introduction}

In 1928 Dirac discovered a matrix valued first-order differential operator as a square root of the Laplacian operator in order to understand elementary particles. Since then this operator and its
various analogues are called Dirac operators in the scientific community. They have played an important role in
physics as well as in mathematics; as an example, let us mention the famous Atiyah-Singer Index Theorem.
Applications to representation theory of reductive Lie groups were started by Parthasarathy \cite{P}, who defined a
Dirac operator on the symmetric space $G/K$ and used it for geometric realization of 
most discrete series representations. This approach was further developed by Atiyah and Schmid \cite{AS}, who constructed 
all discrete series representations as kernels of the Dirac operator acting on the associated spinor bundles.

On the other hand, Lie algebra cohomology (with respect to nilpotent subalgebras) has been a fundamental invariant of
representations of Lie algebras and groups starting from Kostant's 1962 paper \cite{K1}, which gives an algebraic proof of 
the Bott-Borel-Weil theorem.
Lie algebra cohomology is related to characters, asymptotics of matrix coefficients, embeddings into various ``standard"
modules and various geometric realizations of representations. 

It has been clear for decades that the Dirac operators are formally similar to the differentials of the de Rham or Dolbeault cohomology.  There have been two problems with using that analogy in representation theory.  First, the Dirac operators on non-symmetric homogeneous spaces are not as nicely behaved.  Second, the index of the Dirac operator behaves well only for some unitary representations, like the discrete series; so the algebraic tools of representation theory do not work well with it.  The first problem was resolved by Kostant \cite{Kcubic}, who introduced the modified cubic Dirac operator that works well also on non-symmetric homogeneous spaces.  The second problem was resolved by Vogan \cite{V2}, who introduced the concept of Dirac cohomology which applies both to unitary and nonunitary representations: it is the kernel of the
Dirac operator $D$ modulo the intersection of the image and the kernel. He made a conjecture on the infinitesimal 
character of irreducible representations with nonzero Dirac cohomology, which was proved in \cite{HP1}. Kostant generalized this result to the case of his cubic Dirac operator and applied it to the topology of homogeneous spaces \cite{Kdircoh}.
Recently, Kumar \cite{Ku} and Alekseev and Meinreken \cite{AM} found further generalizations
of the results in \cite{HP1} and \cite{Kdircoh}. They put these results into the broader setting 
of non-commutative equivariant cohomology.  

The aim of this paper is to explore the relation between  Dirac cohomology and (nilpotent) Lie algebra cohomology, in those situations when both can be defined. 
Our results show that Dirac cohomology coincides with the corresponding Lie algebra cohomology for a large family of unitary representations, including the discrete series representations, and for all unitary representations in certain special cases.
On the other hand, Dirac cohomology seems to have a better algebraic behavior.   

We have also been able to use Dirac cohomology
and the proved Vogan's conjecture to obtain improvements of some classical results.
In \cite{HP2} we describe how to simplify certain parts of Atiyah-Schmid's construction of 
discrete series representations \cite{AS} and sharpen the Langlands' formula on automorphic forms \cite{L}, \cite{HoP}. The details will appear in a forthcoming book by the first and second named
authors. This suggests that Dirac cohomology may be useful also in other situations, as a new 
tool for tackling problems not accessible to classical cohomology theories.

Let us now describe our main results on Dirac cohomology and its relation to (nilpotent) 
Lie algebra cohomology  more precisely.
Let $G$ be a connected semisimple Lie group with finite center and complexified Lie algebra
$\frg$. Let $\theta$ be a Cartan 
involution of $G$ (and $\frg$), let $K=G^\theta$ be the corresponding
maximal compact subgroup, and let $\frg=\frk\oplus\frp$ be the
complexified Cartan decomposition. As usual, the corresponding real
forms
will be denoted by $\frg_0$, $\frk_0$, etc. (One could also work in a more general setting of a reductive group $G$ in the Harish-Chandra class.)

  If $G$ is of hermitian type, then $\frp$ decomposes as a sum of two abelian subalgebras,
$\frp^+$ and $\frp^{-}$. We prove a Hodge decomposition for $\frp^{-}$-cohomology and 
$\frp^{+}$-homology for all unitary $(\frg,K)$-modules, and show that they are both 
isomorphic to the Dirac cohomology up to a twist by a modular character.
We note that in this setting Enright \cite{E} found an explicit formula for the 
$\frp^{-}$-cohomology of unitary highest weight modules.

For general $\frg$, and any two reductive subalgebras $\frr_1\subset\frr_2\subset\frg$ 
to which the Killing 
form restricts nondegenerately, we show that Kostant's cubic Dirac operator $D(\frg,\frr_1)$ can be decomposed as a sum of two anti-commuting Dirac operators $D(\frg,\frr_2)+D_\Delta(\frr_2,\frr_1)$.
Here $\Delta$ denotes a certain ``diagonal embedding". Decompositions of this kind also appear in \cite{AM}.

In particular, we consider the case when $\frr\subset\frk$ is a
reductive
subalgebra, which is real, i.e., $\frr$ is the complexification of 
$\frr_0=\frr\cap\frg_0$.
Then we can show that for any admissible $(\frg,K)$-module, the Dirac cohomology with respect to 
$D(\frg,\frr)$ is the same as the kernel of $D_\Delta(\frk,\frr)$ on the Dirac cohomology with 
respect to $D(\frg,\frk)$. 
Furthermore, the Dirac cohomology with respect to $D(\frg,\frr)$ is
also the same as the Dirac cohomology with respect to $D(\frg,\frk)$
of the kernel of $D_\Delta(\frk,\frr)$. 

In particular, if $\frg$ and $\frk$ have equal rank, then   
$\frr$ can be a Levi subalgebra $\frl$ of 
a $\theta$-stable parabolic subalgebra $\frqq=\frl\oplus\fru$ of $\frg$, 
with $\frl\subset\frk$.
It follows that the Dirac cohomology with respect to $D(\frg,\frl)$, which is closely related to 
the $\bar\fru$-cohomology, has the advantage of being much easier to calculate. For instance, 
let $\frt$ be a compact Cartan subalgebra, and $\frt\oplus\frn$ be a Borel subalgebra. Then 
we show how to (easily) explicitly calculate the Dirac cohomology of the discrete series 
representations with respect to the Dirac operator $D(\frg,\frt)$. Comparing the obtained
result with Schmid's formula for $\bar\frn$-cohomology in \cite{S}, we see that they are the
same up to an expected modular twist.
In case $G$ is of hermitian type, the Dirac cohomology with respect to $D(\frg,\frl)$ coincides with the $\bar \fru$-cohomology or $\fru$-homology (up to a modular twist) for all unitary representations.

We now describe the organization of our paper. In Sections 2 and 3 we decompose the Kostant cubic Dirac operator as sum of two ``half Dirac operators", which correspond to the differentials of 
$\bar\fru$-cohomology and $\fru$-homology.  Section 3 also contains a new proof of the
Casselman-Osborne theorem on Lie algebra cohomology using an approach similar to that of 
\cite{HP1}. This is
aimed at explaining the formal similarity of the two results.
In Section 4 we prove a Hodge decomposition for $\frp^{-}$-cohomology or $\frp^+$-homology for unitary representations in the Hermitian case, with the Dirac cohomology providing the
``harmonic representatives" of both. The same proof applies to any
finite-dimensional representation of an arbitrary semisimple group. In
Section 5 we define the relative Dirac operators and show that in some
cases Dirac cohomology can be calculated in stages. In Section 6 we develop
this further and show how to use it to calculate
Dirac cohomology explicitly.
In particular, we compare the Dirac cohomology for a Levi subalgebra 
with the $\bar\fru$-cohomology in some cases. 
In Section 7 we obtain a Hodge decomposition and equality of Dirac and $\bar\fru$-cohomology
for arbitrary unitary modules in the Hermitian case.
We conclude the paper by showing that the homological properties of Dirac cohomology are quite
different from those of Lie algebra cohomology. In fact, under certain conditions, we show that
there is a six-term exact sequence of Dirac cohomology corresponding to a short exact sequence
of $(\frg,K)$-modules. So Dirac cohomology resembles a K-theory, rather than a cohomology
theory. 

This work was initiated by David Vogan \cite{V3}. We would like to thank him for many interesting
and stimulating conversations. 
We believe that the results in this paper are not the end of the theory of Dirac cohomology in representation theory, but rather the beginning of further investigations
and applications. For example, the results of this paper should be related to the results
of Connes and Moscovici \cite{CM} in a similar way as the results of \cite{HP1} are related
to \cite{AS}.

\section{Construction of certain Dirac operators}

Let $\frg$ be a complex semisimple Lie algebra, $\frqq=\frl\oplus\fru$  a 
parabolic subalgebra, $\bar\frqq=\frl\oplus\bar\fru$ the opposite parabolic 
subalgebra, and $\frs=\fru\oplus\bar\fru$. Then
$$
\frg=\frl\oplus \frs.
$$
Furthermore, the restrictions of the Killing form $B$ 
to $\frl$ and $\frs$ are non-degenerate, and the above decomposition
is orthogonal. Since 
$\fru$ and $\bar\fru$ are isotropic subspaces in perfect duality
under $B$, we can identify $\bar\fru$ with $\fru^*$; this
identification is $\frl$-invariant. Let $u_1,\dots, u_n$
be a basis of $\fru$, and let $u^*_1,\dots,u^*_n$ be the dual basis of
$\bar\fru$. 

Let $C(\frs)$ be the Clifford algebra of $\frs$. Unlike in \cite{HP1}, we
will use the same defining relations as Kostant, namely 
$$
vw+wv = 2B(v,w);
$$
in particular, if $B(v,v)=1$, then $v^2=1$, and not $-1$ like in  \cite{HP1}.
Of course, over $\bbC$ there is no substantial difference between the
two conventions.

We are going to make use of the well known principle of constructing
invariants by contracting dual indices. The aim is to construct a family
of interesting $\frl$-invariants in $U(\frg)\otimes C(\frs)$. These will
include Kostant's cubic Dirac operator $D$, but we will get $D$ as a sum
of four members of the family, and we will also be able to combine them
in different ways, to get other operators with properties similar to the
properties of $D$. For example, we will have nice expressions for their 
squares. The form of this principle we need is the statement of the following
lemma; the proof is quite easy and essentially reduces to the fact that under
the identification $\hom(\fru,\fru)\cong \fru^*\otimes \fru$, the identity map 
corresponds to the sum $\sum_i u_i^*\otimes u_i$. 

\begin{lemma} Let
$$
\psi:\,\frs^{\otimes 2k} \rightarrow U(\frg)\otimes C(\frs)
$$ 
be a linear map which is $\frl$-equivariant with respect to the adjoint actions. Then  
$$
\sum_I \, \psi(u_I\otimes u^*_I) \in U(\frg)\otimes C(\frs)
$$ 
is independent of the chosen basis $u_i$ and $\frl$-invariant. Here $I=(i_1,\dots,i_k)$
ranges over all $k$-tuples of integers in $\{1,\dots,n\}$, 
$u_I=u_{i_1}\otimes\dots\otimes u_{i_k}$,
and $u^*_I=u^*_{i_1}\otimes\dots\otimes u^*_{i_k}$.
\qed\end{lemma}

For example, $\psi$ can be composed of the obvious inclusions $\frs\hookrightarrow \frg\hookrightarrow U(\frg)$ and $\frs\hookrightarrow C(\frs)$, products, commutators in
$\frg$ and the Killing form $B(.,.)$.  
Here are several examples of this kind which we will study in the following:

\begin{examples} 
\begin{align*}
& A =\sum_i u^*_i\otimes u_i; \\
& A^- =\sum_i u_i\otimes u^*_i; \\
& 1\otimes a = -\frac{1}{4} \sum_{i,j} 1\otimes [u^*_i,u^*_j]u_iu_j; \\
& 1\otimes a^- =-\frac{1}{4} \sum_{i,j} 1\otimes [u_i,u_j]u^*_iu^*_j; \\
& E=1\otimes e = -\frac{1}{2} \sum_i 1\otimes u^*_iu_i.
\end{align*}
\end{examples}

Note the symmetry obtained by exchanging the roles of $\fru$ and $\bar\fru$.
To see how the Dirac operator fits in here, note that one can build an
orthonormal basis $(Z_i)$ of $\frs$ from $u_i$ and $u^*_i$, by putting
\[
Z_j:= \frac{u_j+u^*_j}{\sqrt 2}, \quad Z_{n+j}:=
\frac{i(u_j-u^*_j)}{\sqrt 2},
\]
for $j=1,\dots,n$. Then it is easy to check that
$$
A+A^- = \sum_{i=1}^{2n} Z_i\otimes Z_i.
$$
Also, we can rewrite $a$, $a^-$ and $e$ as follows:

\begin{lemma}
\begin{align*}
& a=-\frac{1}{2} \sum_{i<j}\sum_k B([u^*_i,u^*_j],u_k)\, u_i\wedge u_j\wedge 
u^*_k=-\frac{1}{4} \sum_{i,j} u_iu_j[u^*_i,u^*_j]  ;\\
& a^- = -\frac{1}{2} \sum_{i<j}\sum_k B([u_i,u_j],u^*_k)\, u^*_i\wedge u^*_j
\wedge  u_k = -\frac{1}{4} \sum_{i,j}u^*_iu^*_j [u_i,u_j]; \\
& e = -\frac{1}{2} \sum_i(-u_i u^*_i+2) = \frac{1}{2} \sum_i u_i u^*_i - n.
\end{align*}
\end{lemma}

\begin{proof} Since $[u^*_i,u^*_j]\in\bar \fru$, we can write it as

\[ \sum_k B([u^*_i,u^*_j], u_k)\, u^*_k.\]

\no Also, the sum in the definition of $a$ is clearly twice the same sum over
only those $i,j$ for which $i<j$. 
The only issue is thus to pass from the Clifford product to the wedge product.
For this, we use (1.6) in \cite{Kcubic}. First, since $\fru$ is isotropic, 
$u_iu_j=u_i\wedge u_j$. Next, we calculate

\[ u^*_k (u_i\wedge u_j) = u^*_k\wedge u_i\wedge u_j + B(u^*_k,u_i)\, u_j -
B(u^*_k,u_j)\, u_i.\]

\no The second two terms here are clearly zero if $k$ is different from $i$ and
$j$. For $k=i$, the second term is $u_j$ while the third is zero. 
For $k=j$, the second term is zero while the third is $-u_i$.
It follows that we will be done if we can show

\[  \sum_{i,j}\left( B([u^*_i,u^*_j],u_i)\, u_j -B([u^*_i,u^*_j],u_j)\,  u_i \right)=0
.\]

However, using Lemma 2.1, we see that
$\sum_{i,j} B([u^*_i,u^*_j],u_i)\, u_j$ is an $\frl$-invariant element of
$\fru$. Since there are no nonzero
$\frl$-invariants in $\fru$, this sum must be $0$. Analogously,
$\sum_{i,j} B([u^*_i,u^*_j],u_j)\, u_i =0.$

So we proved the first equality for $a$. Now in the form with wedge product,
we can clearly commute $u_k^*$ in front of $u_i$ and $u_j$, and then we
obtain the second equality by reversing the above argument.

The formulas for $e$ are obvious from the defining relations of $C(\frs)$. 
\end{proof}

Consider now the basis $(b_j)_{j=1,\ldots,2n}$ of $\frs$, given by

$$b_1=u_1,\ldots ,b_n=u_n,\,  b_{n+1}=u^*_1,\ldots ,b_{2n}= u^*_n;$$

\no the dual basis is then

$$d_1=u^*_1,\ldots ,d_n=u^*_n,\,  d_{n+1}=u_1,\ldots ,d_{2n}= u_n.$$

\no Notice that for any $i,j,k$

\[ B([u_i,u_j],u_k) =  B([u^*_i,u^*_j],u^*_k)=0. \]

\no Hence we can write Kostant's cubic element $v$ as

\[ v=-\frac{1}{2}\sum_{1\leq i<j<k\leq 2n} B([d_i,d_j],d_k) \; b_i\wedge b_j\wedge
b_k = a+a^-. \]

\no In particular, we obtain Kostant's cubic Dirac operator as

\[ D= A+A^-+1\otimes(a+a^-). \]

We are also particularly interested in the elements
\[ 
C=A+1\otimes a; \quad C^-=A^-+1\otimes a^-; \quad\text{and}\quad D^-
=C-C^-.
\] 

\no Note that $D=C+C^-$.
We will use the fact that commuting with $E$ operates on
$C,C^-,D$ and $D^-$ in the following way:

\begin{prop} \[ [E,C]=C;\quad [E,C^-]=-C^-;\quad [E,D]=D^- \quad\text{
and}\quad [E,D^-]=D. \]
\end{prop}

\begin{proof} The second two relations follow from the first two, and the first
two are immediate from the following lemma.
\end{proof}

\begin{lemma} Commuting with $e$ in the Clifford algebra $C(\frs)$ 
acts as $I$ on $\fru$ and as $-I$ on $\bar\fru$.
\end{lemma}

\begin{proof} Clearly, for $j\neq i$, $[u^*_j u_j,u_i]=0$. For $j=i$, we
calculate

\[ u^*_i u_i^2 - u_iu^*_i u_i = -u_i(-u_iu^*_i + 2) =-2 u_i. \]

Namely, since $\fru$ is isotropic, $u^2=0$ for any $u\in \fru$.
The first claim now follows, and the second is analogous.
\end{proof}

Kostant \cite{Kcubic}, Theorem 2.16, has calculated
\[
D^2 = \Omega_\frg\otimes 1 - \Omega_{\frl_\Delta} +C.
\]
Here $\Omega_\frg$ denotes the Casimir element of $Z(\frg)\subset U(\frg)$. Further,  
$\Omega_{\frl_\Delta}$ is the Casimir element for the diagonal copy $\frl_\Delta$ of
$\frl$, embedded into $U(\frg)\otimes C(\frs)$ via
\[
X\longmapsto X\otimes 1+1\otimes\alpha(X),\qquad X\in\frl,
\]
where $\alpha:\frl\to\frs\fro(\frs)\to C(\frs)$ is the action map followed by the standard
inclusion of $\frs\fro(\frs)$ into $C(\frs)$ using the identification $\frs\fro(\frs)\cong\bigwedge^2\frs$. Finally, $C$ is the constant $||\rho||^2-||\rho_\frl||^2$.

 Using this result and the above remarks, we
can now quickly calculate the squares of $C$, $C^-$ and $D^-$:

\begin{prop} \[ C^2=(C^-)^2=0 \quad\text{and}\quad (D^-)^2=-D^2. \]
\end{prop}

\begin{proof} From Kostant's expression for $D^2$, it is clear that $D^2$
commutes with all $\frl$-invariant elements of $U(\frg)\otimes C(\frs)$.
In particular, $D^2$ commutes with $E$, and using Proposition 2.4 we see

\[ DD^-+D^-D = [E,D^2]=0. \]

Since $D+D^-=2C$ and $D-D^-=2C^-$, it follows that

\[ 4C^2 = (D+D^-)^2 = D^2+(D^-)^2 = (D-D^-)^2 =4(C^-)^2, \]

so $C^2$=$(C^-)^2$. On the other hand, Proposition 2.4 implies that

\[ [E,C^2]=2C^2 \quad\text{and}\quad [E,(C^-)^2] =-2(C^-)^2. \]

So we see that $C^2=(C^-)^2=0$, and also $D^2+(D^-)^2=4C^2=0$.
\end{proof}

Now we can completely describe the Lie superalgebra $\caD$ spanned by $E$, $C$,
$C^-$ and $\Omega=D^2$ inside the superalgebra $U(\frg)\otimes C(\frs)$.
Here $U(\frg)\otimes C(\frs)$ is an associative superalgebra  with
$\bbZ_2$-grading of the Clifford factor; so it is also a Lie superalgebra
in the usual way, with the supercommutator $[a,b]=ab-(-1)^{\deg a\deg b} ba$.

$\caD$ is a subalgebra of the Lie superalgebra $U(\frg)\otimes C(\frs)$, with the 
commutation relations from
Propositions 2.4 and 2.6. Namely, $\Omega$ is central, and 

\[ [E,C]=C,\quad [E,C^-]=-C^-,\quad [C,C]=[C^-,C^-]=0
\quad\text{and}\quad [C,C^-]=\Omega. \]

The last relation is obtained as follows: 
$\Omega=D^2=(C+C^-)^2=CC^-+C^-C=[C,C^-]$.
Note that we can regard $\caD$
as a $\bbZ$-graded Lie superalgebra with, $C^-,E,\Omega$ and $C$ of
degrees $-1$, $0$, $0$ and $1$ respectively. Note also that the
subalgebra of $\caD$ spanned by $C$, $C^-$ and $\Omega$ is a Heisenberg superalgebra.

The Lie superalgebra $\caD$ can be identified with $\frg\frl(1,1)$, i.e., the
endomorphisms of the superspace $\bbC\oplus\bbC$, with the first $\bbC$ even and 
the second $\bbC$ odd. In Kac's classification \cite{Kac} it is denoted by $\frl(1,1)$. 
It was extensively used by physicists under the name supersymmetric algebra. 
It is a completely solvable Lie superalgebra, and its irreducible finite-dimensional
representations are described in \cite{Kac}. 

To finish this section, let us note that $D$ is independent not
only of the choice of basis $(u_i)$ but also of the choice of 
 $\fru\subset\frs$. On the other hand, $E$, $C$ and $C^-$ do depend on the
choice of $\fru$.

\section{$\bar\fru$-cohomology and $\fru$-homology}

We retain the notation from previous sections. Let $V$ be an admissible $(\frg,K)$-module.
For $p\in \bbN$, let

\[ C^p(\bar \fru,V):= \Hom ({\textstyle\bigwedge^p} \bar \fru,V),  \]

\no be the set of $p$-cochains of the complex defining the $\bar \fru$-cohomology 
of $V$. The differential $d:C^p(\bar \fru,V)  \rightarrow C^{p+1}(\bar \fru,V)$
is given by the usual formula 

\begin{multline*}
 (d\omega)(X_0\wedge \ldots \wedge X_p)=  \sum_{i=0}^p
(-1)^i\;  X_i \cdot \omega (X_0\wedge \ldots 
\wedge \hat{X_i}\wedge \ldots  \wedge X_p) +\\
\sum_{0\leq i<j\leq p}(-1)^{i+j} \; \omega([X_i,X_j] \wedge X_0\wedge \ldots 
\wedge \hat{X_i}\wedge \ldots \wedge \hat{X_j}\wedge \ldots \wedge X_p)
\end{multline*}

We have the following identifications :

\[   C^p(\bar \fru,V) =  \Hom ({\textstyle\bigwedge^p} \bar \fru,V)\simeq 
\Hom ({\textstyle\bigwedge^p} (\fru^*),V) 
\simeq   \Hom (({\textstyle\bigwedge^p} \fru)^*,V)\simeq
V \otimes  {\textstyle\bigwedge^p} \fru 
\] 

Like in the previous section, we fix a basis $(u_i)_{i=1,\ldots,n}$ of $\fru$ and
denote the dual basis of $\bar\fru$ by $(u^*_i)_{i=1,\ldots,n}$. 

\begin{lemma}
Through the above identifications, the differential $d:V\otimes{\textstyle\bigwedge^p} \fru  \rightarrow V \otimes  {\textstyle\bigwedge^{p+1}} \fru $ is given by 
\begin{multline*}
d(v\otimes Y_1\wedge \ldots  \wedge Y_p)=  \sum_{i=1}^n   
u^*_i\cdot v\otimes u_i\wedge Y_1 \wedge \ldots \wedge Y_p\\
 +\frac{1}{2}\sum_{i=1}^n \sum_{j=1}^p v\otimes  u_i\wedge Y_1 \wedge \ldots 
\wedge [u^*_i,Y_j]_\fru \wedge \ldots \wedge Y_p 
\end{multline*}

\no where $[u^*_i,Y_j]_\fru$ denotes the projection of $[u^*_i,Y_j]$ on 
$\fru$.
\end{lemma}

\begin{proof}
This is a straigtforward calculation, starting from the fact that the identification
$\bigwedge^p (\fru^*)=(\bigwedge^p \fru)^*$ is given via 
$(f_1\wedge\dots\wedge f_p)(X_1\wedge\dots\wedge X_p) = \det f_i(X_j)$.
\end{proof}

The space $V\otimes \bigwedge^p \fru $ is also the space of $p$-chains for the 
$\fru$-homology of $V$, with differential $\partial: V\otimes {\textstyle\bigwedge^p} \fru  \rightarrow V\otimes {\textstyle\bigwedge^{p-1}} \fru$ 
given by :
\begin{multline*}
 \partial(v\otimes Y_1\wedge \ldots  \wedge Y_p )=  \sum_{i=1}^p (-1)^i   
Y_i\cdot v \otimes Y_1 \wedge \ldots \hat Y_i \wedge \ldots
\wedge Y_p  + \\
 \sum_{1\leq i<j\leq p}  (-1)^{i+j} v\otimes [Y_i,Y_j]\wedge Y_1 \wedge \ldots 
\wedge \hat Y_i \wedge \ldots  \wedge \hat Y_j \wedge \ldots Y_p 
\end{multline*}

Note that we are tensoring with $\bigwedge^\cdot\fru$ from the right and not from the
left as usual; this is because we will have an action of $U(\frg)\otimes C(\frs)$
on $V\otimes\bigwedge^\cdot\fru$ which will be more natural in this order.

For reasons that will become apparent later, we will instead of
$\partial$ consider the operator $\delta=-2\partial$.
Of course, $\delta$ defines the same homology as $\partial$.

To get our Dirac operators act, we need to consider the
$U(\frg)\otimes C(\frs)$ - module $V\otimes S$, where $S$ is the spin
module for the Clifford algebra $C(\frs)$.
We will use the identification of $\bigwedge^\cdot\fru$ with $S$,
given explicitly in \cite{Kdircoh} and \cite{Kgen}. Namely, one can construct $S$ as 
the left ideal
in $C(\frs)$ generated by the element $u^*_{\top} = u^*_1\dots u^*_n$.
One then has $S=(\bigwedge ^\cdot\fru)u^*_{\top}$, which is isomorphic to
$\bigwedge^\cdot\fru$ as a vector space, and the action of
$C(\frs)$ is given by left Clifford multiplication. Explicitly, 
$u\in\fru$ and $u^*\in\bar\fru$ act on 
$Y=Y_1\wedge\dots\wedge Y_p\in\bigwedge^p\fru$ by

\[ u\cdot Y = u\wedge Y\]

\[ u^*\cdot Y = 2 \sum_{i=1}^p(-1)^{i+1} B(u^*,
Y_i)Y_1\wedge\dots\hat Y_i\dots\wedge Y_p. \] 

Namely, since $\fru$ and $\bar\fru$ are isotropic, the Clifford and
wedge products coincide on each of them; in particular, $u^*u^*_{\top}=0$.

The natural action of $\frl$ on $V\otimes S$ is the tensor product of
the restriction of the $\frg$-action on $V$ and the spin action on
$S$.
On the other hand, the usual $\frl$ action on $\bar\fru$-cohomology and
$\fru$-homology is given by the adjoint action on $\bigwedge^\cdot\bar\fru$
and $\bigwedge ^\cdot\fru$. Thus, our identification of $\bigwedge
^\cdot\fru\otimes V$ with $V\otimes S$ is not an
$\frl$-isomorphism. However, as was proved in \cite{Kgen}, Prop.3.6,
the two actions differ only by a
twist with the one dimensional $\frl$-module $Z_{\rho(\bar\fru)}$ 
of weight $\rho(\bar\fru)$. 

This means that, if we consider $d$ and $\delta$ as operators on
$V\otimes S$ via the above identification, then as an $\frl$-module,
the cohomology of $d$
gets identified with $H^\cdot(\bar\fru,V) \otimes Z_{\rho(\bar\fru)}$, while
the homology of $\delta$ gets identified with $H_\cdot(\fru,V) \otimes
Z_{\rho(\bar\fru)}$.

\begin{prop} Under the action of $U(\frg)\otimes C(\frs)$ on $V\otimes S$, 
the operators $C=A+1\otimes a$ and $C^-=A ^-+1\otimes a ^-$ 
from Section 2 act as $d$ and
$\delta$ respectively. In particular, the cubic Dirac operator $D=C+C^-$
acts as $d+\delta = d-2\partial$.
\end{prop}

\begin{proof} We are going to use the explicit formulas for the action of
$C(\frs)$ on $S$ given above. From these formulas, it is immediate that the
action of $A$ coincides with the first (single) sum in the expression
for $d$, while the action of $A ^-$ transforms $v\otimes Y_1\wedge\dots\wedge
Y_p \in V\otimes \bigwedge ^p \fru$ into

\[ 
\sum_{i=1}^n  u_i v\otimes 2\sum_{k=1}^p (-1)^{k+1}
B(u_i^*,Y_k)\, Y_1\wedge\dots\hat Y_k\dots \wedge Y_p. 
\]

Since $\sum_iB(u^*_i,Y_k)\,  u_i = Y_k$, we see that this is equal to
minus twice the first (single) sum in the expression for the $\fru$-homology
operator $\partial$.

It remains to identify the action of the cubic terms $a$ and $a ^-$. 

We use the expression for $a$ from Lemma 2.3, i.e., $a=-\frac{1}
{4}\sum_{i,j}u_i u_j [u^*_i,u^*_j]$. This element transforms
$v\otimes Y_1\wedge\dots\wedge
Y_p\in V\otimes \bigwedge ^p \fru$ into

\begin{align*}
-\frac{1}{4} v\otimes \sum_{i,j} u_iu_j\; 2\sum_{k=1}^p (-1)^{k+1} B([u_i^*,u^*_j],Y_k)\, 
Y_1\wedge\dots\hat Y_k\dots\wedge Y_p   \\
=\frac{1}{2} v\otimes \sum_{i,j,k} (-1)^{k+1}B([u^*_i,Y_k],u^*_j)\, u_i\wedge
u_j\wedge  Y_1\wedge\dots\hat Y_k\dots\wedge Y_p.
\end{align*}

Now we sum $\sum_j B([u^*_i,Y_k],u^*_j)\, u_j =[u_i^*,Y_k]_\fru$, 
and after commuting  $[u^*_i,Y_k]_\fru$ into its proper place,
we get the second (double) sum in the expression for $d$.

For $a ^-$ we use its definition from 2.2. To write the calculation
nicely, we introduce the following notation: for
$Y=Y_1\wedge\dots\wedge Y_p$, let 

\[  
\hat Y_{k,l} = Y_1\wedge\dots\hat Y_k\dots\hat Y_l\dots \wedge Y_p, \qquad
\text{if}\; k<l.
\]

If $k>l$ then we change the sign and define 

\[
\hat Y_{k,l} =- Y_1\wedge\dots\hat Y_l\dots\hat Y_k\dots\wedge Y_p, \qquad
\text{if}\; k>l.
\]

If $k=l$, we set $\hat Y_{k,l}=0$. This now allows us to write

\[
u^*_iu^*_j \cdot Y_1\wedge\dots\wedge Y_p  =4 \sum_{k,l}(-1)^{k+l}
B(u^*_i,Y_k)B(u^*_j,Y_l) \hat Y_{k,l}.
\]

It follows that $1\otimes a ^-$ transforms $v\otimes Y_1\wedge\dots\wedge
Y_p$ into

\[
-\frac{1}{4}\;4\;v\otimes \sum_{i,j,k,l} (-1)^{k+l}
B(u^*_i,Y_k)B(u^*_j,Y_l)\, [u_i,u_j]\wedge \hat Y_{k,l}.
\]

Upon summing up $\sum_i B(u^*_i,Y_k)\, u_i = Y_k$ and 
$\sum_j B(u^*_j,Y_l)u_j = Y_l$, we get that this is equal to

\[
-v\otimes \sum_{k,l} (-1)^{k+l}[Y_k,Y_l]\wedge \hat Y_{k,l}.
\]

This is now clearly invariant for exchanging the roles of $k$ and $l$,
hence it is twice the same sum extending just over $k<l$, i.e., minus
twice the second (double) sum in the expression for $\partial$.
\end{proof}

It is now clear why we considered $\delta=-2\partial$ instead of just
$\partial$; in this way we have the action of $D$ being equal to
$d+\delta$.

Before we go on, let us note how the
element $E$ of Section 2 acts on $V\otimes S$; it is in fact a degree
operator up to a shift.
This means $E$ can be used 
to identify the degrees in which the cohomology (homology) is
appearing.

\begin{prop} The element $E$ of $U(\frg)\otimes C(\frs)$ 
from 2.2 acts on $V\otimes \bigwedge
^k \fru $ as multiplication by the scalar $k-n$. Consequently,
it preserves the kernel and image of both $d$ and $\delta$, and
hence acts on the $k$-th cohomology of $d$ and the $k$-th homology of 
$\delta$,  by the same scalar $k-n$. In particular, $E+n$ is the degree operator.
\end{prop}

\begin{proof} Using Lemma 2.5 we see that for any $u\in\fru$,
$eu=ue+u$.
On the other hand, by Lemma 2.3, $e=\frac{1}{2}\sum_i u_iu^*_i
-n$, hence $eu^*_{\top} = -nu^*_{\top}$. It now immediately follows that in
$C(\frs)$
we have $e(Y_1\wedge\dots\wedge
Y_k\;u^*_{\top})=(k-n)Y_1\wedge\dots\wedge Y_k\;u^*_{\top}$, so the
action of $E$ on $V\otimes \bigwedge^k \fru$ is indeed 
multiplication by the scalar $k-n$.

It now immediately follows that
$E$ preserves the kernel and image of $d$ and $\delta$, as these are
homogeneous operators (of degree $1$ and $-1$ respectively).
(Note that this last assertion can also be obtained from the
commutation relations of Proposition 2.4.)
\end{proof}

We will now state a result for the operators $C$ and $C^-$
analogous to the one obtained
 for $D$ in \cite{HP1} and \cite{Kdircoh}. A corollary will be the Casselman-Osborne Theorem.
Our goal here is not to give a new proof of the Casselman-Osborne Theorem, the existing ones being
completely satisfactory, but to shed some light on the formal similarity between
the Casselman-Osborne Theorem and the main result of \cite{HP1}.

Define $ d_D, d_C, d_{C^-}:\; 
U(\frg)\otimes  C(\frs)\rightarrow U(\frg)\otimes  C(\frs) $ by 

\begin{align*}
 d_D(x)&=Dx-\epsilon_x xD\\
 d_C(x)&=Cx-\epsilon_x xC\\
 d_{C^-}(x)&=C^-x-\epsilon_x xC^-
\end{align*}

\no where $\epsilon_x$ is 1 for even $x$ and $-1$ for odd $x$. In the following we fix a compact
group $L$ with complexified Lie algebra $\frl$.

\begin{thm}
$ d_D, d_C$ and $ d_{C^-}$ are $L$-equivariant. They induce maps from 
$ ( U(\frg)\otimes  C(\frs))^L$ into itself
and $d_D^2=d_C^2=d_{C^-}^2=0$ on 
$ ( U(\frg)\otimes  C(\frs))^L$.

Furthermore $Z(\frl_\Delta)\subset \ker  d_D$ (resp.  
$\ker  d_C$, $\ker  d_{C^-}$), and one has
$\ker  d_D= Z(\frl_\Delta)\oplus\im d_D$, (resp.  
$\ker  d_C= Z(\frl_\Delta)\oplus\im d_C$,
 $\ker  d_{C^-}= Z(\frl_\Delta)\oplus\Im d_{C^-}$).
\end{thm}

\begin{proof} The result for $ d_D$ is due to Kostant \cite{Kdircoh}.
The proof is the same as the
proof of the main result of \cite{HP1}. We give details for $d_C$, the proof
for $ d_{C^-}$ being entirely similar. As in \cite{HP1}, we use the standard  
filtration on $ U(\frg)$, 
which induces a filtration  $(F_nA)_{n\in \bbN}$ on 
$A:=  U(\frg)\otimes  C(\frs)$.
 This filtration being $L$-invariant, 
it induces in turn a filtration on $A^L$. Clearly, 
\[ C=\sum_{j=1}^n  u_j^* \otimes u_j+ (1\otimes a) \in F_1 A^L. \]

The $\bbZ_2$-gradation on the Clifford algebra $ C(\frs)$ induces a $\bbZ_2$-gradation on $A$. 
We set $A^0=  U(\frg)\otimes  C(\frs)^0$ and $A^1=  U(\frg)\otimes  C(\frs)^1$.
 Then $A=A^0\bigoplus A^1$ and 
this $\bbZ_2$-gradation is compatible with the filtration   $(F_nA)_{n\in \bbN}$

If $x\in F_n A^0$, then $ d_C(x)= d_C^0(x)=Cx-xC \in F_{n+1}A^1$. If $x \in F_n A^1$, then 
$ d_C(x)= d_C^1(x)=Cx+xC \in F_{n+1}A^0$. Thus 

\[  d_C^0:\;  F_n A^0 \rightarrow  F_{n+1}A^1 \quad \text{ and }  \quad d_C^1:\;  F_n A^1 \rightarrow  F_{n+1}A^0 \]

\no induce

\[\bar   d_C^0:\;  \Gr_n \, A^0 \rightarrow   \Gr_{n+1}A^1 \quad \text{ and }  
\quad \bar  d_C^1:\;    \ \Gr_n\, A^1 \rightarrow  \Gr_{n+1}\, A^0. \]

Let $\bar C= \sum_{j=1}^n  u_j^* \otimes u_j \in S^1(\frg)\otimes C(\frs)$ be the image of 
$C$ in  $\left(\Gr_1\, A^1\right)^L$
(notice that the cubic term disappears since it is in $F_0 A^1$).
If $x \in  F_n A^0$, 

\[ \bar   d_C^0 (\bar x)=\overline{Cx-xC}=\bar C \bar x-\bar x \bar C\]

\no and if  $x \in  F_n A^1$, 

\[ \bar   d_C^1 (\bar x)=\overline{Cx+xC}=\bar C \bar x+\bar x \bar C.\]

\no Therefore $\bar d_C= \bar  d_C^0 \oplus \bar  d_C^1: \Gr A \rightarrow \Gr A$. Note also that unlike $d_D$, $d_C$ is a differential on all of $U(\frg)\otimes C(\frs)$,
because $C^2=0$. Hence $\bar d_C$ is a differential on $S(\frg)\otimes C(\frs)$.

Let us compute $\bar d_C(\bar x)$ for 
$\bar x= y \otimes u_{i_1}\cdots  u_{i_k}u_{j_1}^*\cdots u_{j_l}^* \in S(\frg)\otimes  C(\frs)$.
We can assume $i_r$ are different from each other and likewise for $j_s$.

\begin{multline*}
\bar  d_C(\bar x)=\left(\sum_{j=1}^n  u_j^* \otimes u_j \right) \left( y \otimes u_{i_1}\cdots  u_{i_k}
u_{j_1}^*\cdots  u_{j_l}^* \right) - \\
(-1)^{k+l} \left( y \otimes u_{i_1}\cdots  u_{i_k}
 u_{j_1}^*\cdots u_{j_l}^* \right)\left(\sum_{j=1}^n  u_j^* \otimes u_j \right) \\
= \sum_{j=1}^n  u_j^* y \otimes \left( (-1)^k  u_{i_1}\cdots  u_{i_k} u_j
 u_{j_1}^*\cdots  u_{j_l}^*-(-1)^{k+l} u_{i_1}\cdots  u_{i_k}
 u_{j_1}^*\cdots  u_{j_l}^*u_j \right)
\end{multline*}

If $j\neq j_s$ for all $s$, then the contribution to the sum is zero. If $j=j_s$, then

\begin{multline*}
u_{j_s} u_{j_1}^*\cdots  u_{j_l}^*-(-1)^l u_{j_1}^*\cdots  u_{j_l}^*u_{j_s} = 
(-1)^{s-1} u_{j_1}^*\cdots  u_{j_{s-1}}^* (u_{j_s} u_{j_s}^*+u_{j_s}^* u_{j_s})u_{j_{s+1}}^*\cdots  u_{j_l}^* \\ = 2(-1)^{s-1} u_{j_1}^*\cdots \widehat{u_{j_s}^*} \cdots  u_{j_l}^* . 
\end{multline*}

So we see 

\[\bar  d_C(\bar x)= \sum_{s=1}^l -2(-1)^{k+s} u_{j_s}^* y\otimes  u_{i_1}\cdots  
u_{i_k} u_{j_1}^*\cdots \widehat{u_{j_s}^*}\cdots u_{j_l}^*.\]

Since $\frg=\frl\oplus \fru \oplus \bar \fru$ and $ C(\frs)\simeq \bigwedge^\cdot \frs =
\bigwedge^\cdot \fru  \otimes 
 \bigwedge^\cdot\bar \fru$,
one has

\[ S(\frg)\otimes   C(\frs)= S(\frl)\otimes S(\fru)\otimes {\textstyle\bigwedge^\cdot} \fru   
\otimes  S(\bar \fru)  \otimes {\textstyle\bigwedge^\cdot} \bar \fru. \]

It follows that $\bar  d_C$ is (up to a sign depending on the 
$\bigwedge^\cdot\fru$-degree)
equal to 
$-2\Id \otimes d_{\bar \fru}$, where $d_{\bar \fru}$ is the Koszul
differential for the vector space $\bar \fru$. It is well known that $d_{\bar \fru}$ is exact except at degree zero,
where the cohomology is $\bbC$, embedded as the constants. It follows that  $\bar d_C$  is exact except at
 degree zero, where the cohomology is  $S(\frl)\otimes S(\fru)\otimes 
\bigwedge^\cdot \fru$ embedded in the obvious way.
It remains to pass to the invariants:

\begin{lemma}
The differential $\bar d_C$ on  $\left( S(\frg)\otimes   C(\frs) \right)^L$  is exact except 
at degree zero. The zeroth cohomology is  
$\left( S(\frl)\otimes S(\fru)\otimes\bigwedge^\cdot \fru \right)^L=S(\frl)^L \otimes 1$, 
embedded in the obvious way. More precisely

\[\ker\bar  d_C= \left( S(\frl)^L \otimes 1\right)\oplus \im 
\bar  d_C \]

\end{lemma}

To prove the lemma, we need to 
show that $\left(S(\frl)\otimes S(\fru)\otimes \bigwedge^\cdot \fru \right)^L=
 S(\frl)^L \otimes 1\otimes 1$. To see this, we may choose some element $h$ in the
 center of $\frl$, such that $\frl$ is the centralizer of $h$ in $\frg$, $\ad h$ has 
real eigenvalues, and $\fru$ is the sum of the eigenspaces of $\ad h$ corresponding to
the positive eigenvalues. 
Making $h$ act on an element in $\left( S(\frl)\otimes 
S(\fru)\otimes \bigwedge^\cdot \fru \right)^L$, we see that 
this element has to be in $S(\frl)^L \otimes 1\otimes 1$.

We can now finish the proof of the theorem. We proceed as in \cite{HP1} by induction on the degree of the filtration.
\end{proof}

 Let $z\in Z(\frg)$. Since $z\otimes 1\in \ker  d_C$, we can write

\[  z\otimes 1=\eta_\frl(z)+Ca+aC   \]

\no for some 
$a\in \left(  U(\frg) \otimes   C(\frs) \right)^L$ and $\eta_\frl(z)\in Z(\frl_\Delta)$.

Our goal is now to compute $\eta_\frl(z)$. Let $\frh$ be a Cartan subalgebra of $\frl$, and let us denote respectively by $W_\frl$ and $W_\frg$ the Weyl groups of $\frh$ in $\frl$ and $\frg$.

We have Harish-Chandra isomophisms 

$$\mu_{\frl/\frh}:\, Z(\frl)\rightarrow S(\frh)^{W_\frl},\qquad
\mu_{\frg/\frh}:\, Z(\frg)\rightarrow S(\frh)^{W_\frg},$$

\no and an obvious inclusion $i:\, S(\frh)^{W_\frg}\rightarrow S(\frh)^{W_\frl}$.
Set $\mu_{\frg/\frl}:= \mu_{\frl/\frh}^{-1}\circ i \circ \mu_{\frg/\frh}  $.
With this notation we have:

\begin{lemma}
For all $z \in Z(\frg)$, $\eta_{\frl}(z)=\mu_{\frg/\frl}(z)$.
\end{lemma}

\begin{proof}  The proof is similar to the proof of Theorem 4.2 in \cite{Kdircoh}, 
but much simpler. We give only a sketch. Let $\frb$ be a Borel subalgebra of $\frg$ containing $\bar \fru$, and  suppose the Cartan algebra $\frh$ has been chosen to lie in $\frb$.

Let $V_\lambda$ be the irreducible finite dimensional  representation with highest weight
$\lambda$ (relative to $\frb$), and let $v_\lambda$ be a non-zero highest weight vector
in $V_\lambda$. Recall the element $u^*_{\top}$ in $C(\frs)$ used to define the spin module $S$.

One can see easily that $C\cdot(v_\lambda \otimes u^*_{\top})=0$ and  
that $v_\lambda \otimes u^*_{\top} \in V_\lambda \otimes S$ defines a non-zero cohomology class
in $H^0(\bar \fru,V_\lambda)$. Since the infinitesimal character of $V_\lambda$  is 
given by $\lambda+\rho$, and any $q\in Z(\frl_\Delta)$ acts by the scalar
$\mu_{\frl/\frh}(q)(\lambda+\rho)$ on $v_\lambda\otimes u^*_{\top}$ (see \cite{Kdircoh}, 
Theorem 4.1), we get 

\[  \mu_{\frg/\frh}(z)(\lambda+\rho)=  \mu_{\frl/\frh}(\eta_\frl(z))(\lambda+\rho) \]
for all $z \in Z(\frg)$.

Since this is true for all dominant weights $\lambda$, and since these form a Zariski dense
set in $\frh^*$, we conclude that indeed  $\mu_{\frg/\frh}(z)=\mu_{\frl/\frh}(\eta_\frl(z))$.
\end{proof}

Let $V$ be a representation of $\frg$ admitting an infinitesimal character $\chi_V$ and 
let $z\in  Z(\frg)$. Then $z$ acts on $H^\cdot(\bar \fru,V)$ by the scalar $\chi_V(z)$.

Notice that $z\otimes 1$ acts on  $\ker C/\im C$ as $\eta_\frl(z)$. Namely, $Ca+aC$
leaves $\ker C$ and $\im C$ stable, and induces the zero action on $\ker C/\im C$.

Thus, the induced action of $\eta_\frl(z)$ on $\ker C/\im C\cong H^\cdot(\bar \fru,V)\otimes Z_{\rho(\bar\fru)}$ is equal to the scalar multiplication by $\chi_V(z)$. This is exactly the 
statement of the Casselman-Osborne Theorem (see \cite{CO}, or \cite{V1}, Theorem 3.1.5).
Namely, our definition of $\mu_{\frg/\frl}=\eta_\frl$ differs from the map $Z(\frg)\to Z(\frl)$ 
from the Casselman-Osborne Theorem exactly by the above $\rho$-shift.

\section{Hodge decomposition for $\frp^-$ - cohomology}

In this section we study positive definite hermitian forms on $V\otimes S$ for unitary modules $V$,
such that $d$ and $\delta$ are minus adjoint to each other. Here unitarity is defined with respect
to a fixed real form $\frg_0$ of $\frg$; we also fix a corresponding Cartan involution $\theta$ and
Cartan decompositions $\frg_0=\frk_0\oplus\frp_0$ and $\frg=\frk\oplus\frp$.

Such a form however exists only in very special situations; if $\frg$ is simple, then the pair 
$(\frg,\frk)$ must be hermitian symmetric, $\frl$ must be equal to $\frk$, and $\fru$ and $\bar\fru$ 
must be the abelian subalgebras $\frp^\pm$ of $\frp$. In that case, we get a Hodge type decomposition:
$\frp^+$-homology, $\frp^-$-cohomology and the Dirac cohomology are all represented by the space of
``harmonics", $\ker D=\ker D^2$.

A parallel result holds much more generally if $V$ is a finite-dimensional module. Since this case
is known, we include it only as a remark. However we start by studying a sufficiently general situation
which applies to both cases.

Let $\frr$ be a subalgebra of $\frg$ to which $B$ restricts non-degenerately 
and let $\frs$ be the orthogonal complement to $\frr$ with respect to
$B$. We choose a maximal isotropic subspace $\frs^+$ of $\frs$. Since 
$$
\langle X,Y\rangle =-2B(X,\theta \bar Y)
$$
(with $\,\,\bar{ }\,\,$ denoting conjugation with respect to $\frg_0$)
defines a positive definite Hermitian form on $\frg$ and hence also on $\frs^+$, 
the subspace $\frs^-=\theta\overline{\frs^+}$ intersects $\frs^+$ trivially. Let 
$S=\bigwedge^\cdot\frs^+$ be a spin module for the Clifford algebra $C(\frs)$ 
corresponding to this polarization. We extend the form $\langle\,,\rangle$ to 
all of $S$ in the usual way, using the determinant. 

\begin{lemma} With respect to the form $\langle X,Y\rangle$ on $S$, 
the adjoint of the operator $X\in\frs\subset C(\frs)$ is $-\theta\bar X$. 
\end{lemma}

\begin{proof} We need to show that
$$
\langle X\cdot\lambda,\mu\rangle = \langle\lambda,-\theta\bar X\cdot\mu\rangle
$$
for any $X\in\frs$ and $\lambda,\mu\in S$. Assume $X\in\frs^+$. We can assume that 
$\lambda=\lambda_1\wedge\dots\wedge\lambda_s$ and $\mu=\mu_1\wedge\dots\wedge\mu_t$,
with $X$, $\lambda_i$ and $\mu_j$ all belonging to a fixed orthonormal basis of $\frs^+$.
Then both sides of the equality to be established are nonzero precisely when $X$ is
one of the $\mu_j$'s, while the other $\mu_j$'s are precisely (all) the $\lambda_i$'s.
In this case both sides are clearly the same.

If $X\in\frs^-$ the claim follows by considering $\theta\bar X\in\frs^+$. In case $\frs$
is odd-dimensional, the extra element is handled easily.
\end{proof}

We now assume that $\frr$ is real and $\theta$-stable, i.e., $\frr$ is the complexification 
of a subalgebra $\frr_0\subset\frg_0$ with $\frr_0=\frr_0\cap\frk_0\oplus\frr_0\cap\frp_0$.

\begin{lemma} (a) The cubic part $v$ of Kostant's cubic Dirac operator $D$ is self-adjoint
with respect to the form $\langle\,,\rangle$ on $S$.

(b) If $\frr=\frl$ for a parabolic subalgebra $\frl\oplus\fru$, so that $v=a+a^-$ as in the
previous sections, then $a$ is adjoint to $a^-$.
\end{lemma}

\begin{proof} (a) Let us choose bases $Z_i$ of $\frs_0\cap\frk_0$ and $Z_r'$ of $\frs_0\cap\frp_0$
orthonormal with respect to $\langle\,,\rangle$. By Lemma 4.1, the
adjoint of $Z_i$ is $-Z_i$ and the adjoint of $Z'_r$ is $Z_r'$.
Moreover, the dual bases of $Z_i$, $Z_r'$
with respect to $B$ are $-Z_i$, $Z_r'$, and so 
\begin{multline*}
v=-\frac{1}{2}\left(\sum_{i<j<k} B([-Z_i,-Z_j],-Z_k)Z_iZ_jZ_k
  +\sum_{i;r<s} B([-Z_i,Z_r'],Z_s')Z_iZ_r'Z_s'\right)  \\
= \frac{1}{2}\left(\sum_{i<j<k} B([Z_i,Z_j],Z_k)Z_iZ_jZ_k
  +\sum_{i;r<s} B([Z_i,Z_r'],Z_s')Z_iZ_r'Z_s' \right)  
\end{multline*}
(the other terms are zero). Since the adjoint of $Z_iZ_jZ_k$ is
$(-Z_k)(-Z_j)(-Z_i)=Z_iZ_jZ_k$, the adjoint of 
$Z_iZ_r'Z_s'$ is $Z_r'Z_s'(-Z_i)=Z_iZ_r'Z_s'$, and the coefficients
are real, we see that $v$ is self-adjoint. 

(b) By the formulas from Lemma 2.3, it is enough to show that 
$u_iu_j[u_i^*,u_j^*]$ is adjoint to $[u_i,u_j]u_i^*u_j^*$. This
however follows immediately from Lemma 4.1: we can normalize the
$u_i$'s so that $u_i^*=-\theta \bar u_i$, hence the adjoint of 
$[u_i,u_j]u_i^*u_j^*$ is 
$u_ju_i(-\theta\overline{[u_i,u_j]})=u_iu_j[u_i^*,u_j^*]$
\end{proof}

Let now $V$ be a unitary $(\frg,K)$-module and let $\frqq=\frl\oplus\fru$ be
a $\theta$-stable parabolic subalgebra. Then the adjoint of the operator
$u_i$ on $V$ with respect to the unitary form is $-\bar u_i$, while the
adjoint of $u_i$ on $S$ is either $-\bar u_i$, or $\bar u_i$, depending on
whether $u_i$ is in $\frk$ or in $\frp$. If
we consider the tensor product hermitian form on $V\otimes S$,
denoted again by $\langle\,,\rangle$, and as before 
denote by $A_\frk$ and $A_\frp$ the $\frk$ - respectively $\frp$ - parts of $A$, then we see:

\begin{cor} With respect to the form $\langle\,,\rangle$ on
$V\otimes S$, the adjoint of $A_\frk$ is $A_\frk^-$ while the adjoint of $A_\frp$ is $-A_\frp^-$.
Hence the adjoint of $C=A_\frk+A_\frp+1\otimes a$ is $A_\frk^-+1\otimes a^- - A_\frp^-$
and the adjoint of
$C^-=A_\frk^-+A_\frp^-+1\otimes a^-$ is $A_\frk+1\otimes a -A_\frp$. \qed
\end{cor}

We will use this corollary in Sections 5 and 7. Now we turn our attention to the case when
$A_\frk$ and $1\otimes a$ do not appear in $C$. This is the already mentioned case,
when $\frl=\frk$, and $\fru$ is contained in $\frp$. Then $\fru$ is forced to be abelian
and we denote as usual $\fru=\frp^+$, $\bar\fru=\frp^-$. In this case, the Dirac operator
$D=D(\frg,\frl)$ is the ``ordinary" Dirac operator corresponding to $\frk$, there is no
cubic part, and we conclude

\begin{cor} Let $(\frg,\frk)$ be a hermitian symmetric pair and set $\frl=\frk$. Let $V$ be
a unitary $(\frg,K)$-module and consider the form $\langle \,,\rangle$ on $V\otimes S$.
Then the operators $C=d$ and $C^-=\delta$ are minus adjoints of each other. Hence the 
Dirac operator $D=D(\frg,\frl)=D(\frg,\frk)$ is anti-self-adjoint.

More generally, if $\frl$ contains $\frk$ (this can happen when $\frg$ is not simple),
then $D(\frg,\frl)$ is anti-self-adjoint with respect to $\langle \,,\rangle$. \qed
\end{cor}

For the rest of this section we assume that $\frl=\frk$.
So $D$ is anti-self-adjoint with respect to the  
positive definite form $\langle\,,\rangle$ on $V\otimes S$. 
In particular, the
operators $D$ and $D^2$ have the same kernel on $V\otimes S$.

By \cite{P}, Proposition 3.2 (or more generally, \cite{Kcubic}, Theorem 2.16), we know that 
\[
D^2=\Omega_\frg\otimes 1 - \Omega_{\frk_\Delta} +C,
\]
where $\Omega_\frg$ and $\Omega_{\frk_\Delta}$ are the Casimir operators for
$\frg$ and diagonally embedded $\frk$, and $C$ is the constant 
$||\rho_\frg||^2-||\rho_\frk||^2$. 
It follows that if $\Omega_\frg$ acts on $V$ by a constant, 
then $\Omega_{\frk_\Delta}$ is up to a
constant equal to $D^2$ on $V\otimes S$. Since $\Omega_{\frk_\Delta}$ acts by a scalar on
each $\tilde K$-type in $V\otimes S$, the same is true for $D^2$. So $D^2$
is a semisimple operator, i.e., $V\otimes S$ is a direct sum of eigenspaces for $D^2$.
In particular:

\begin{cor} If the $(\frg,K)$-module $V$ has infinitesimal character,
then $V\otimes S= \ker(D^2)\oplus \im(D^2)$
\end{cor}

\begin{proof}  We have seen that $V\otimes S$ is a direct sum of
eigenspaces for $D^2$. Clearly, the zero eigenspace is $\ker D^2$, and
the sum of the nonzero eigenspaces is $\im D^2$.
\end{proof}

It is now easy to obtain a variant of the usual Hodge decomposition.
The following arguments are well known; see e.g. \cite{W}, Scholium 9.4.4.
We first have

\begin{lemma} (a) $\ker D = \ker d\cap\ker\delta$;

(b)   $\im\delta$ is
orthogonal to $\ker d$ and $\im d$ is orthogonal to $\ker\delta$.  
\end{lemma}

\begin{proof}
(a) Since $D=d+\delta$, it is clear that $\ker
d\cap\ker\delta$ is contained in $\ker D$. On the other hand, if
$Dx=0$, then $dx=-\delta x$, hence $\delta dx =-\delta^2x=0$. So
$\langle dx,dx\rangle=\langle -\delta dx,x\rangle=0$,
hence $dx=0$. Now $Dx=0$ implies that also $\delta x=0$.

(b) is obvious since $d$ and $\delta$ are minus adjoint to each other.
\end{proof}

Combining Corollary 4.5, Lemma 4.6 and the fact $\ker D=\ker D^2$, we get

\begin{thm} Let $(\frg,\frk)$ be a hermitian symmetric pair and set $\frl=\frk$ and $\fru=\frp^+$.
Let $V$ be an irreducible unitary $(\frg,K)$-module. Then:

(a) $V\otimes S = \ker D\oplus \im d\oplus \im\delta$;

(b) $\ker d = \ker D\oplus \im d$;

(c) $\ker\delta=\ker D\oplus \im\delta$. 

In particular, the Dirac cohomology of $V$ is equal to $\frp^-$-cohomology and to
$\frp^+$-homology, up to modular twists:
\[
\ker D\cong H^\cdot(\frp^-,V)\otimes Z_{\rho(\frp^-)}\cong H_\cdot(\frp^+,V)\otimes Z_{\rho(\frp^-)}.
\]
 More precisely, (up to modular twists) the Dirac cohomology $\ker D$
is the space of harmonic representatives for both $\frp^-$-cohomology and 
$\frp^+$-homology. \qed
\end{thm}

\begin{rmk}
There is a variant of the above results for a finite dimensional module $V$.
In this case one puts a positive
definite form on $V$ invariant under the compact form $\frk_0\oplus i\frp_0$
of $\frg$; this is sometimes called
an admissible form. 
In this way the adjoint of $u_i$ on $V$ will be $\bar u_i$ if
$u_i\in\frp$ and $-\bar u_i$ if $u_i\in \frk$. Hence the adjoint of
$u_i$ on $V$ is $u_i^*$ for all $i$. 

So this form combines with $\langle\,,\rangle$
on $S$ better than a unitary form. In particular,  
$d$ and $\delta$ are now adjoint to each other with respect to $\langle\,,\rangle$, 
$D$ is self-adjoint, and the above arguments, including Theorem 4.7, 
apply without change. There is no need here to assume that $\frl=\frk$; 
$\frqq=\frl\oplus\fru$ can be any parabolic subalgebra of (any) $\frg$, with the
assumption that $\frl$ is $\theta$-stable and real.

This case was however already known; it is implicit in 
\cite{Kgen} and it is explicitly mentioned in \cite{V3}. 

Moreover, we can generalize this further and prove self-adjointness of $D$ for the 
more general setting when $\frr\subset\frg$ is not necessarily a Levi subalgebra of a 
parabolic subalgebra,
but any real and $\theta$-stable reductive subalgebra (to which $B$ then
restricts nondegenerately). Namely, getting back to the setting of Lemma 4.2 (a)
and its proof, we see that not only the cubic term $v$ but also the 
linear term $\sum Z_i\otimes Z_i +\sum Z_r'\otimes Z_r'$ of $D$
is self adjoint with respect to $\langle\,,\rangle$ on $V\otimes S$.
Hence $D$ is self-adjoint.
\end{rmk}

\begin{rmk} Here is an example which shows that a Hodge decomposition like in Theorem
4.7 does not hold generally if $\frl$ is not compact, even for unitary modules in the 
Hermitian symmetric case.
Let $W_0$ be the $(\frg,K)$-module of the even Weil representation of $Sp(4,\bbR)$, and 
let $\frl$ be the $\frg\frl(2)$ containing the compact Cartan subalgebra and 
corresponding to the short noncompact root. Then $\frl$ is a Levi factor of a 
$\theta$-stable parabolic subalgebra with abelian $\fru$. One can explicitly calculate 
everything and show that the Dirac cohomology is strictly larger than $\fru$-homology
or $\bar\fru$-cohomology (while the latter two are still equal).
\end{rmk}

\section{Relative Dirac operators}

In this section we compare various Dirac operators arising from
two compatible decompositions
$$
\frg=\frr\oplus\frs = \frr'\oplus\frs'.
$$
Here both decompositions are like in the Kostant's setting for the cubic Dirac operator:
$\frr$ and $\frr'$ are reductive subalgebras of $\frg$ such that the Killing form 
$B$ restricts to a nondegenerate form on each of them, while $\frs$ respectively
$\frs'$ are the orthogonals of $\frr$ respectively $\frr'$.
Compatibility of the two decompositions means
$$
\frr=\frr\cap\frr'\oplus \frr\cap\frs';\qquad \frs=\frs\cap\frr'\oplus \frs\cap\frs';\qquad
\frr'=\frr'\cap\frr\oplus \frr'\cap\frs;\qquad \frs'=\frs'\cap\frr\oplus \frs'\cap\frs.
$$
Clearly, all these decompositions are orthogonal.   

\begin{example} If $\frl$ is a $\theta$-stable Levi subalgebra of $\frg$, then $\frr=\frl$
and $\frr'=\frk$ satisfy the above conditions.
\end{example}

Now $\frr\cap\frr'$ is also a reductive subalgebra of $\frg$ to which $B$ restricts
nondegenerately. The corresponding decomposition is
$$
\frg=\frr\cap\frr' \oplus (\frs\oplus \frr\cap\frs').
$$
The corresponding Dirac operator will be denoted by $D(\frg,\frr\cap\frr')$.
We note that this setting includes the case in which $\frr'$ is a subalgebra of $\frr$,
and therefore the setting $\frr_1 \subset \frr_2 \subset \frg$ from the introduction.

To write down $D(\frg,\frr\cap\frr')$,
take orthonormal bases $Z_i$ for $\frs$, and $Z'_j$ for $\frr\cap\frs'$. Identify

\begin{equation}
U(\frg)\otimes C(\frs\oplus \frr\cap\frs') = U(\frg)\otimes C(\frs)\bar\otimes C(\frr\cap\frs'),
\end{equation}

\no where $\bar\otimes$ denotes the $\bbZ_2$-graded tensor product. Now if $W_k$ is the union of
the bases $Z_i$ and $Z'_j$, then by Kostant's definition, $D(\frg,\frr\cap\frr')$ is the element
$$
D(\frg,\frr\cap\frr') = \sum_k W_k\otimes W_k -\frac{1}{2}\sum_{i<j<k} B([W_i,W_j],W_k)\otimes 
W_i W_j W_k 
$$
of $U(\frg)\otimes C(\frs\oplus \frr\cap\frs')$.

Kostant's original definition uses exterior multiplication
instead of Clifford multiplication in the second sum (with the Clifford and 
exterior algebras identified as vector spaces via
Chevalley identification).
For orthogonal vectors, there is however no difference between exterior and Clifford
multiplication, so the above definition is the same as Kostant's.

Taking into account the decomposition (1), we see

\begin{multline}
D(\frg,\frr\cap\frr') = \sum_i Z_i\otimes Z_i\otimes 1 + \sum_j Z_j'\otimes 1\otimes Z_j' 
-\frac{1}{2}
\sum_{i<j<k} B([Z_i,Z_j],Z_k)\otimes Z_i Z_j Z_k\otimes 1 - \\
\frac{1}{2} \sum_{i<j}\sum_k  B([Z_i,Z_j],Z_k')\otimes Z_i Z_j\otimes Z_k' -\frac{1}{2}
\sum_{i<j<k} B([Z_i',Z_j'],Z_k')\otimes 1\otimes Z_i' Z_j' Z_k'. 
\end{multline}

Note that the terms with $Z_i$, $Z_j'$ and $Z_k'$ do not appear, because
$B([Z_i,Z_j'],Z_k')=B(Z_i,[Z_j',Z'_k])=0$, as $[Z'_j,Z'_k]\in\frr$ is orthogonal to $\frs$.

We consider $U(\frg)\otimes C(\frs)$ as the subalgebra $U(\frg)\otimes C(\frs)\otimes 1$
of $U(\frg)\otimes C(\frs)\bar\otimes C(\frr\cap\frs')$. In view of this, we see that
the first and third sum in (2) combine to give $D(\frg,\frr)$, the Kostant's cubic Dirac 
operator corresponding to $\frr\subset\frg$.

The remaining three sums come from the cubic Dirac operator
corresponding to $\frr\cap\frr'\subset\frr$. However, this is an element
of the algebra $U(\frr)\otimes C(\frr\cap\frs')$, and this algebra has to be embedded
into $U(\frg)\otimes C(\frs)\bar\otimes C(\frr\cap\frs')$ diagonally, not as
$U(\frr)\otimes 1\otimes C(\frr\cap\frs')$. Namely, we use the diagonal embedding
$U(\frr)\cong U(\frr_\Delta)\subset U(\frg)\otimes C(\frs)$ from the setting
$\frg=\frr\oplus\frs$; the definition is the same as for $\frr=\frl$ above Proposition 2.6.
Thus we embed
$$
\Delta: U(\frr)\otimes C(\frr\cap\frs') \cong U(\frr_\Delta)\bar\otimes C(\frr\cap\frs')
\subset U(\frg)\otimes C(\frs)\bar\otimes C(\frr\cap\frs').
$$

We will denote $\Delta(D(\frr,\frr\cap\frr'))$ by $D_\Delta(\frr,\frr\cap\frr')$ and call it
a relative Dirac operator. In case when there are several subalgebras and confusion might arise,
we will use more precise notation $\Delta_{\frg,\frr}$ instead of $\Delta$ and give up the
notation $D_\Delta(.)$.
This diagonal embedding has already been
used by Kostant and Alekseev-Meinrenken; in particular, decompositions like
the following one can be founded in \cite{AM}.

\begin{thm} With notation as above,

(i) $D(\frg,\frr\cap\frr') = D(\frg,\frr) + D_\Delta(\frr,\frr\cap\frr')$;

(ii) The summands $D(\frg,\frr)$ and $D_\Delta(\frr,\frr\cap\frr')$ anticommute.
\end{thm}

\begin{proof} To prove (i), we need to describe the image under $\Delta$ of

\begin{equation}
D(\frr,\frr\cap\frr')= \sum_i Z_i'\otimes Z_i' -\frac{1}{2}\sum_{i<j<k} B([Z_i',Z_j'],Z_k')\otimes Z_i'Z_j'Z_k' \in U(\frr)\otimes C(\frr\cap\frs'),
\end{equation}

\no and see that it matches the second, fourth and fifth sum in (2). In fact, it is obvious that
the image under $\Delta$ of the second sum in (3) equals the fifth sum in (2), and it remains to 
identify

\begin{equation}
\sum_i \Delta(Z_i'\otimes Z_i') = \sum_i Z_i'\otimes 1\otimes Z_i'+ 1\otimes\alpha(Z_i')\otimes Z_i'.
\end{equation}

\no Namely, $\Delta(Z\otimes Z')=Z\otimes 1\otimes Z'+1\otimes\alpha(Z)\otimes Z'$, where 
$\alpha: \frr\to \frs\fro(\frs)\hra C(\frs)$ is the action map of $\frr$ on $\frs$, followed
by the standard inclusion of $\frs\fro(\frs)$ into $C(\frs)$, given by
$$
E_{ij}-E_{ji} \mapsto \frac{1}{2} Z_iZ_j,
$$
where $E_{ij}$ are the matrix units relative to the basis $Z_i$.

Thus we are left with showing that the second sum in (4) equals the third sum in (2), i.e., that
$$
\sum_k 1\otimes\alpha(Z_k')\otimes Z_k' =-\frac{1}{2} \sum_{i<j}\sum_k  B([Z_i,Z_j],Z_k')\otimes Z_i Z_j\otimes Z_k'.
$$
This will follow if we see that
$$
\alpha(Z_k') = -\frac{1}{2} \sum_{i<j}  B([Z_i,Z_j],Z_k') Z_i Z_j
$$
for any $k$. But this last equality is clear from the definition of $\alpha$. Namely, the
matrix of $\ad Z_k'$ on $\frs$ in the basis $Z_i$ is 
$\sum_{i<j} B([Z_k',Z_j],Z_i)(E_{ij}-E_{ji})$.

To prove (ii), we use the fact that $D(\frg,\frr)$ commutes with $\frr_\Delta$, which is
one of the most basic properties of $D(\frg,\frr)$. It follows that the anticommutator
$$
\big[D(\frg,\frr)\otimes 1,(Z_i'\otimes 1+1\otimes\alpha(Z_i'))\otimes Z_i'\big] =
\big[D(\frg,\frr), Z_i'\otimes 1+1\otimes\alpha(Z_i')\big]\otimes Z_i' 
$$
is zero for any $i$. Hence $\big[D(\frg,\frr)\otimes 1,\Delta(\sum_i Z_i'\otimes Z_i')\big]=0$.
It remains to see that also
$$
\big[D(\frg,\frr)\otimes 1, 1\otimes1\otimes (-\frac{1}{2})\sum_{i<j<k}  B([Z_i',Z_j'],Z_k')Z_i'Z_j'Z_k'\big]=0.
$$
This follows from the definition of $\bar\otimes$, since all the $C(\frs)$-parts of
the monomial terms of $D(\frg,\frr)$, and also all $Z_i'Z_j'Z_k'\in C(\frr\cap\frs')$, are odd.
\end{proof}

\begin{example}
In the setting of Example 5.1, we obtain
$$
D(\frg,\frl\cap\frk) = D(\frg,\frl)+ D_\Delta(\frl,\frl\cap\frk) = D(\frg,\frk)+D_\Delta(\frk,\frl\cap\frk),
$$
and both decompositions have anticommuting summands. If $\frl$ is contained in $\frk$, which is
possible if and only if $\frg$ and $\frk$ have equal ranks, then we get only one nontrivial
decomposition, $D(\frg,\frl)=D(\frg,\frk) +D_\Delta(\frk,\frl)$. This case will be of special
interest below.
\end{example}

We now want to use Theorem 5.2 to relate the Dirac cohomology of the various Dirac operators
involved. In some cases one can apply the following easy fact from linear algebra.
We define the cohomology of any linear operator $T$ on a vector space $V$ to be the vector space
$H(T)=\ker T/(\im T\cap\ker T)$.

\begin{lemma} Let $A$ and $B$ be anticommuting linear operators on a vector space $V$.
Assume that $A^2$ diagonalizes on $V$, i.e. $V=\bigoplus_\lambda V_\lambda$, with $A^2=\lambda$
on $V_\lambda$. Then the cohomology $H(A+B)$ of $A+B$ on $V$ is the same as the cohomology of the
restriction of $A+B$ to $V_0=\ker A^2$.
\end{lemma}

\begin{proof} Since $A+B$ commutes with $A^2$, its kernel, image and cohomology decompose
accordingly to eigenspaces $V_\lambda$. We thus have to prove that $A+B$ has no cohomology
on $V_\lambda$ for $\lambda\neq 0$. In other words, we are to prove that 
$\ker(A+B)\subset\im(A+B)$ on $V_\lambda$.

Let $v\in V_\lambda$ be such that $(A+B)v=0$, i.e., $Av=-Bv$. Then 
$$
(A+B)Av = A^2v+BAv = A^2 v -ABv  = 2A^2 v = 2\lambda v,
$$
and hence $v=\frac{1}{2\lambda}(A+B)Av$ is in the image of $A+B$.
\end{proof}

\begin{cor} In the setting of Lemma 5.4, assume further that $\ker A^2=\ker A=H(A)$; so 
$\ker A\cap\im A=0$. Then $H(A+B)$ is equal to the cohomology of $B$ restricted to the
cohomology (i.e., kernel) of $A$.
\end{cor}

\begin{proof} By Lemma 5.4, $H(A+B)$ is the cohomology of $A+B$ on $\ker A$. But on $\ker A$,
$A+B=B$.
\end{proof}

To apply this to Dirac cohomology, denote by $H_D(\frg,\frr;V)$ the Dirac cohomology of
a $(\frg,K)$-module $V$ with respect to $D(\frg,\frr)$; analogous notation will be used
for other Dirac operators. The reader should bear in mind that $H_D(\frg,\frr;V)$ is in
fact the cohomology of the operator $D(\frg,\frr)$ on the space $V\otimes S$.

\begin{cor} 
Let $\frr\subset\frk$ be a reductive subalgebra of $\frg$, with 
$B|_{\frr\times\frr}$ nondegenerate. As usual, let $\frs$ be the
orthocomplement of $\frr$. Assume that either $\dim\frp$ is even, or
$\dim\frs\cap\frk$ is even\footnote{It should be possible to eliminate this 
assumption by using the graded version of the spin module, as explained at the end of Section 8.}.
 Let $V$ be an irreducible unitary $(\frg,K)$-module. Then  
$$
H_D(\frg,\frr;V) =  H_D(\frk,\frr;H_D(\frg,\frk;V)),
$$
i.e., the Dirac cohomology can be calculated ``in stages", as the $D(\frk,\frr)$-cohomology of
the $D(\frg,\frk)$-cohomology.
\end{cor}

\begin{proof} 
Since $V$ is unitary, we can consider the form $\langle\,,\rangle$ on 
$V\otimes S_\frp$ introduced in Section 4, where $S_\frp$ is the spin
module for $C(\frp)$.
 We can extend this form to all of $V\otimes S$, by choosing
any positive definite form on the spin module $S_{\frs\cap\frk}$ for
$C(\frs\cap\frk)$.
Here we identify $S=S_\frp\otimes S_{\frs\cap\frk}$, which can be
done by the assumption on dimensions.
Let $A=D(\frg,\frk)$ and $B=D_\Delta(\frk,\frr)$.

By Corollary 4.3, $A$ is anti-self-adjoint, and consequently
the conditions of Corollary 5.5 are satisfied. So the cohomology with respect to
$D(\frg,\frr)$ is the cohomology with respect to $B$ of $\ker A = H_D(\frg,\frk;V)\otimes S_{\frs\cap\frk}$.

Now $H_D(\frg,\frk;V)\subset V\otimes S_\frp\subset V\otimes S$ is a $\tilde K$-module, 
with Lie algebra $\frk$ acting through $\frk_\Delta$. The Dirac cohomology of this module
with respect to $D(\frk,\frr)$ is thus identified with the cohomology with respect to
$B=D_\Delta(\frk,\frr)$.
\end{proof}

\section{The case of compact Levi subalgebra}

In this section we first consider a reductive subalgebra $\frr_0$ of
$\frg_0$ contained in $\frk_0$. Later on we will specialize to the
case when $\frr=\frl$ is a Levi subalgebra of a $\theta$-stable parabolic 
subalgebra of $\frg$.

The first thing we will do in this situation is generalize Corollary
5.6 to nonunitary modules. Like there, we assume for simplicity that
either $\dim\frp$ is even or $\dim\frs\cap\frk$ is even, so that we
can write the spin module as $S= S_\frp\otimes S_{\frs\cap\frk}$.
The idea is to reverse the roles of $D(\frg,\frk)$ and $D_\Delta(\frk,\frr)$. Namely, for any 
admissible $(\frg,K)$-module $V$, we can decompose $V\otimes S_\frp$ into a 
direct sum of finite dimensional (unitary) modules for the spin double cover $\tilde K$ of $K$. 
Hence, by Remark 4.8, there is a positive definite  
form $\langle\,,\rangle$ on 
$V\otimes S= V\otimes S_\frp\otimes S_{\frs\cap\frk}$, 
such
that $D_\Delta(\frk,\frr)$ is self-adjoint with respect to $\langle\,,\rangle$.

It follows that $B=D_\Delta(\frk,\frr)$ is a semisimple operator, while for $A=D(\frg,\frk)$ 
we still have that $A^2$ is semisimple. In this situation, we have the following lemma which 
complements Corollary 5.5:

\begin{lemma} Let $A$ and $B$ be anticommuting linear operators on a vector space $V$, such that
$A^2$ and $B$ are semisimple (i.e., can be diagonalized). Then $H(A+B)$ is the cohomology
(i.e., the kernel) of $B$ acting on $H(A)$. 
\end{lemma}

\begin{proof} Applying Lemma 5.4, we can replace $V$ by $\ker A^2$, i.e., assume $A^2=0$.
On the other hand, by Corollary 5.5, $H(A+B)$ is the cohomology of $A$ acting on $\ker B$.

Since $B$ is semisimple, we can decompose 
\[
V=\ker B\oplus\bigoplus_{\lambda} V_{\lambda}\oplus V_{-\lambda}
\]
into the (discrete) sum of eigenspaces for $B$. Here if both $\lambda$ and $-\lambda$ are
eigenvalues, we choose one of them to represent the pair. 
Since $A$ anticommutes with $B$, it preserves $\ker B$, and maps $V_\lambda$ to $V_{-\lambda}$ and vice versa.
Therefore, $H(A)$ decomposes into a $\ker B$-part
and $V_\lambda\oplus V_{-\lambda}$-parts. The $\ker B$-part is equal to $H(A+B)$ and we will be done if we show that $B$
has no kernel on $H(A|_{V_\lambda\oplus V_{-\lambda}})$. 
Let $v=v_1+v_2\in V_\lambda\oplus V_{-\lambda}$ be in $\ker A$, and assume that $Bv\in\im A$.
This implies $\lambda v_1-\lambda v_2$ is in $\im A$, so $v_1-v_2\in\im A$. This however can
only happen if both $v_1$ and $v_2$ are in $\im A$, again because $A$ exchanges $V_\lambda$
and $V_{-\lambda}$. But then also $v=v_1+v_2$ is in $\im A$, so $v$ is zero in cohomology and we are done.

\end{proof}

This now immediately implies the following theorem which says that Dirac cohomology with respect
to a subalgebra $\frr\subset\frk$ as above can be calculated ``in stages".

\begin{thm} Let $\frr_0$ be a reductive
subalgebra of $\frg_0$ contained in $\frk_0$.
Let $V$ be an admissible $(\frg,K)$-module. Then the Dirac cohomology with respect to
 $D(\frg,\frr)$ can be calculated as the Dirac cohomology with respect to $D(\frk,\frr)$ of
the Dirac cohomology with respect to $D(\frg,\frk)$ of $V$. In other words:
\[
H_D(\frg,\frr;V)=H_D(\frk,\frr;H_D(\frg,\frk;V)).
\]
Also, we can reverse the order of taking Dirac cohomology, i.e., 
\[
H_D(\frg,\frr;V)=H(D(\frg,\frk)\big|_{H_D(\frk,\frr;V)}).
\]
\end{thm}

\begin{proof} The first formula was explained above, and the second is a direct application
of Corollary 5.5, with $A=D_\Delta(\frk,\frr)$ and $B=D(\frg,\frk)$ (opposite from Corollary
5.6).
\end{proof}

For the rest of this section we
consider a $\theta$-stable parabolic subalgebra $\frqq=\frl\oplus\fru$
of $\frg$, with the Levi subalgebra $\frl$ contained in $\frk$. In particular, 
there is a Cartan subalgebra $\frt$ of $\frg$ contained in $\frl\subset\frk$; so $\frg$ and
$\frk$ have equal rank. The opposite parabolic subalgebra is $\bar\frqq=\frl\oplus\bar\fru$.
As before, we denote $\frs=\fru\oplus\bar\fru$, so
$\frg=\frl\oplus\frs$.

We apply the above considerations to $\frr=\frl$.
Since $H_D(\frg,\frk;V)$ is a finite dimensional $\tilde K$-module, and $\frk$ and $\frl$
have equal rank, 
$H_D(\frk,\frl;H_D(\frg,\frk;V))$ is given by \cite{Kdircoh}, 
Theorem 5.1. (This can also be read off 
from Remark 4.8, by using Poincar\'e duality to pass to $\fru$-cohomology and then the better
known Kostant's formula for $\fru$-cohomology; in fact, this is how Kostant proves it.)
This gives $H_D(\frg,\frl;V)$ very explicitly
provided we know $H_D(\frg,\frk;V)$ explicitly. 
For example, one can in this way calculate the Dirac cohomology of the discrete
series representations with respect to the (compact) Cartan subalgebra $\frt$:

\begin{example} Let $V=A_\frb(\lambda)$ be a discrete series representation; here 
$\frb=\frt\oplus\frn$ is a Borel subalgebra of $\frg$ containing a compact Cartan subalgebra 
$\frt$. The infinitesimal character of $V$ is $\lambda+\rho$. Then the
Dirac cohomology of $V$ with respect to $D(\frg,\frk)$ consists of a single
$\tilde K$-type $V(\mu)$, whose highest weight is $\mu=\lambda+\rho_n$, where 
$\rho_n=\rho(\fru\cap\frp)$. 
This is obtained from the highest weight of the lowest $K$-type of $V$, $\lambda+2\rho_n$, by
shifting by $-\rho_n$ (the lowest weight of $S$). 

This result is contained in the work of Parthasarathy and Schmid. One can also prove it as follows:
it is shown in \cite{HP1}, Proposition 5.4, that this
$\tilde K$-type is contained in the Dirac cohomology. Since $V$ has a unique lowest $K$-type,
and since $-\rho_n$ is the lowest weight of the spin module, with multiplicity one, it follows
that any other $\tilde K$-type has strictly larger highest weight, and thus can not contribute to 
the Dirac cohomology.

We now apply the above mentioned Kostant's formula (Theorem 5.1 of \cite{Kdircoh}) to calculate 
the Dirac cohomology with respect to $D(\frk,\frt)$ (we again stress that $\frk$ and $\frt$ have
equal rank):
$$
H_D(\frk,\frt;V(\mu))= \Ker D(\frk,\frt)=\bigoplus_{w\in W_\frk}\C_{w(\mu+\rho_\frk)}.
$$
It follows from $\mu+\rho_\frk=\lambda+\rho$ that 
$$H_D(\frk,\frt;V(\mu))=\bigoplus_{w\in W_\frk}\C_{w(\lambda+\rho)}.$$
\end{example}

\begin{rmk} 
Comparing with Schmid's  formula in Theorem 4.1 of \cite{S}, we have
$$
H^*(\bar\frn,A_\frb(\lambda))=H_D(\frg,\frt;A_\frb(\lambda))\otimes \C_{\rho(\bar\frn)}.
$$
(note that Schmid's $\frn$ is our $\bar\frn$, and his $\lambda$ is our $\lambda+\rho$.)

In other words, $\frn$-cohomology of a discrete series representation
coincides with the Dirac cohomology
up to a $\rho$-shift.  This fact is however not covered by our 
results in Sections 4 and 7. This indicates that it should be possible to generalize our results; 
however, the next example shows they do not hold in full generality.
\qed
\end{rmk}

\begin{example} In the setting of Remark 4.9, one can show that 
$H_D(\frl,\frt;H_D(\frg,\frl;W_0))$ is strictly larger than $H_D(\frg,\frt;W_0)$.
Here $\frt$ is the compact Cartan subalgebra, and the last cohomology can be calculated
in stages, as $H_D(\frk,\frt;H_D(\frg,\frk;W_0))$.
\end{example}

We now want to similarly analyze the ``half-Dirac" operators $C$ and $C^-$.
Let $u_1,\dots,u_k$ be a basis for $\fru\cap\frk$ and let $v_1,\dots,v_p$ be a basis for
$\fru\cap\frp$. These can be taken to be the root vectors
corresponding to compact, respectively noncompact positive roots, with respect to some
$\Delta^+(\frg,\frt)$ compatible with $\fru$. We normalize these bases so that
the dual bases for $\bar\fru\cap\frk$ respectively $\bar\fru\cap\frp$ with respect to the Killing form are $u_i^*=-\bar u_i$ respectively $v_i^*=\bar v_i$.

As before, $D=D(\frg,\frl)$, $C=C(\frg,\frl)=A+1\otimes a$ and 
$C^-=C^-(\frg,\frl)=A^-+1\otimes a^-$ denote the Dirac operator for
the pair $(\frg,\frl)$ and its parts. We can now further decompose these parts and write
$$
A=A_\frk+A_\frp;\quad a=a_{\frk}+a_{\frk\frp}+a_\frp
$$
and analogously for $A^-$ and $a^-$. Here           

\begin{gather*}
A_\frk=\sum_i u_i^*\otimes u_i,\quad A_\frp=\sum_i v_i^*\otimes v_i;  \\
a_\frk=-\frac{1}{4}\sum_{i,j} [u_i^*,u_j^*]u_iu_j,\quad 
a_{\frk\frp}=-\frac{1}{2}\sum_{i,j} [u_i^*,v_j^*]u_iv_j,\quad 
a_{\frp}=-\frac{1}{4}\sum_{i,j} [v_i^*,v_j^*]v_iv_j. 
\end{gather*}

The expressions for $A_\frk^-$, $A_\frp^-$, $a_{\frk}^-$, $a_{\frk\frp}^-$ and $a_\frp^-$
are obtained by exchanging $u_i$ with $u_i^*$ and $v_i$ with $v_i^*$. 

In the following, we will consider the Clifford algebra
$C(\frs)$ as a subalgebra of $U(\frg)\otimes C(\frs)$, embedded as $1\otimes C(\frs)$.
In particular, $1\otimes a_\frk$ gets identified with $a_\frk$, $1\otimes a_{\frk\frp}$ 
with $a_{\frk\frp}$ and so on. 

Recall that by Theorem 5.2, $D(\frg,\frl)=D(\frg,\frk)+D_\Delta(\frk,\frl)$, where 
$D_\Delta(\frk,\frl)$ is the image of $D(\frk,\frl)$ under the diagonal embedding
$\Delta: U(\frk)\otimes C(\frs\cap\frk) \to U(\frg)\otimes C(\frp)\bar\otimes C(\frs\cap\frk)$.
Here $\Delta$ sends $1\otimes C(\frs\cap\frk)$ identically onto 
$1\otimes 1\otimes C(\frs\cap\frk)$, and for $X\in\frk$, 
$$
\Delta (X\otimes 1) = X\otimes 1\otimes 1 + 1\otimes\alpha(X)\otimes 1,
$$
with $\alpha:\frk\to C(\frp)$ given by the action of $\frk$ on $\frp$ as before.

Clearly, $D(\frg,\frk)=A_\frp +A_\frp^-$, while $D_\Delta(\frk,\frl)$ is the sum of all other
of the above parts. We want to make this more precise; namely, in the obvious notation,
$D(\frk,\frl)=C(\frk,\frl)+C^-(\frk,\frl)$, and we want to identify the images of these
summands under $\Delta$. Denote these images by $C_\Delta(\frk,\frl)$ and $C^-_\Delta(\frk,\frl)$.

To do this, we use an expression for $\alpha:\frk\to C(\frp)$ in terms of a basis and a dual
basis: if $w_r$ is a basis of $\frp$ with dual basis $w_r^*$, then
$$
\alpha(X)=-\frac{1}{4}\sum_{r,s} B([w_r^*,w_s^*],X) w_rw_s.
$$
(This was already used in Section 5, for an orthonormal basis. The proof of this slightly more
general version is the same.)
Applying this to the basis $v_1,\dots,v_p,v_1^*,\dots,v_p^*$ and the dual basis 
$v_1^*,\dots,v_p^*,v_1,\dots,v_p$, we get

\begin{multline*}
\alpha(X)=-\frac{1}{4}\sum_{j,k}B([v_j^*,v_k^*],X) v_jv_k 
 -\frac{1}{2}\sum_{j,k}B([v_j,v_k^*],X) v_j^*v_k \\
 +\frac{1}{2}\sum_{j}B([v_j,v_j^*],X) 
-\frac{1}{4}\sum_{j,k}B([v_j,v_k],X) v_j^*v_k^* 
\end{multline*}
(we used $v_j^*v_k=-v_kv_j^*+2\delta_{jk}$).

Since $C(\frk,\frl)=\sum_i u_i^*\otimes u_i -\frac{1}{4}\otimes\sum_{i,j}[u_i^*,u_j^*]u_iu_j$, 
we see that
$\Delta(C(\frk,\frl))=A_\frk +\sum_i 1\otimes\alpha(u_i^*)\otimes u_i +a_\frk$.
We need to calculate the middle term, $\sum_i 1\otimes\alpha(u_i^*)\otimes u_i$.
Applying the above expression for $\alpha$, we get four sums over $i,j$ and $k$.

We first notice that the first of these three sums is 0, since $B$ is
0 on $\bar\fru$. To calculate the second sum, write 
$B([v_j,v_k^*],u_i^*)=B(v_j,[v_k^*,u_i^*])$, and observe that since 
$[v_k^*,u_i^*]\in\bar\fru\cap\frp$, 
$\sum_j B(v_j,[v_k^*,u_i^*])v_j^* = [v_k^*,u_i^*]$. Therefore the second of the three sums is
$$
-\frac{1}{2}\otimes\sum_{i,k}[v_k^*,u_i^*] v_ku_i = a_{\frk\frp}.
$$
The third sum is 0, since we can assume $[v_j,v_j^*]$ is in $\frl$ and hence orthogonal to
$u_i^*$. Namely, we can choose $v_j$ and $v_j^*$ (and also $u_i$ and $u_i^*$) to be root vectors
with respect to $\frt$.

Finally, the fourth sum is calculated by noting that since $[v_j,v_k]\in\fru\cap\frk$,
$\sum_i B([v_j,v_k],u_i^*)u_i = [v_j,v_k]$. It follows that the fourth sum is
$$
-\frac{1}{4}\otimes\sum_{j,k} v_j^*v_k^*[v_j,v_k] = a_\frp^-.
$$
A completely analogous calculation applies to $C^-(\frk,\frl)$, so we proved:

\begin{prop} Under the diagonal map 
$\Delta: U(\frk)\otimes C(\frs\cap\frk) \to U(\frg)\otimes C(\frp)\bar\otimes C(\frs\cap\frk)$,
$C(\frk,\frl)$ and $C^-(\frk,\frl)$ correspond to
$$
C_\Delta(\frk,\frl)= A_\frk+a_\frk+a_{\frk\frp} + a_\frp^- \qquad\text{and}\qquad
C^-_\Delta(\frk,\frl)= A^-_\frk+a^-_\frk+a^-_{\frk\frp} + a_\frp. \qed
$$
\end{prop}

Note the unexpected feature of this result, the mixing of the positive and negative
parts under the diagonal embedding. Namely, $a_\frp$ and $a_\frp^-$ have opposite positions
from the ones one would expect. So we do not have an analogue of Theorem 5.2 for $C$ and $C^-$,
unless $a_\frp=a^-_\frp=0$. This last thing happens precisely when the pair $(\frg,\frk)$ is
hermitian symmetric. This is the reason why we are able to obtain results about $\fru$-cohomology
only in hermitian case. Maybe this peculiar behavior has something to do with the fact
that some of the most concrete results
about $\frn$-cohomology, like \cite{E}, \cite{C} or \cite{A}, are also obtained in hermitian
situation only. 

Another difficulty with the non-hermitian case is the fact that while we can write
$D(\frg,\frk)=A_\frp+A^-_\frp$, the two summands here are not differentials (they
are also not $K$-invariant). So there is
no hope to get a Hodge decomposition like the one in Section 4. Of course, there is also
no $\fru\cap\frp$-homology or cohomology, since $\fru\cap\frp$ is not a Lie algebra. Yet there is
perfectly well defined Dirac cohomology for $D(\frg,\frk)$, and one can hope that it
will somehow replace the nonexistent $\frp^-$-cohomology.

What we do get without the hermitian assumption, is a copy of $\frg\frl(1,1)$ inside 
$U(\frg)\otimes C(\frs)$,
spanned by $C_\Delta(\frk,\frl)$, $C_\Delta^-(\frk,\frl)$, 
$E_\frk=\Delta E(\frk,\frl)=-\frac{1}{2}\otimes\sum_i u_i^*u_i$ and $\Delta D(\frk,\frl)^2$. 
This follows immediately from the fact that $\Delta$ is 
(obviously) a morphism of superalgebras. 

If the pair $(\frg,\frk)$ is hermitian symmetric, as we will assume in the following,
then there is another copy of $\frg\frl(1,1)$ inside $U(\frg)\otimes C(\frs)$, supercommuting
with the first one. It is spanned by $A_\frp=C(\frg,\frk)$, $A_\frp^-=C^-(\frg,\frk)$, 
$E_\frp=E(\frg,\frk)= -\frac{1}{2}\otimes\sum_i v_i^*v_i$, and $D(\frg,\frk)^2$.
The fact that these two copies of $\frg\frl(1,1)$ supercommute is completely analogous to
Theorem 5.2.(ii). (Without the hermitian assumption, one could try to replace this second 
$\frg\frl(1,1)$ by a smaller superalgebra, spanned just by $D(\frg,\frk)$ and $D(\frg,\frk)^2$.)

Note that $C(\frg,\frl)=C_\Delta(\frk,\frl)+C(\frg,\frk)$ and 
$C^-(\frg,\frl)=C^-_\Delta(\frk,\frl)+C^-(\frg,\frk)$.
Clearly, an analogous decomposition holds for $E(\frg,\frl)$; it also holds for $D(\frg,\frl)^2$.
Namely, since the two copies of $\frg\frl(1,1)$ supercommute, we have
$$
D(\frg,\frl)^2 = [C(\frg,\frl),C^-(\frg,\frl)] =[C_\Delta(\frk,\frl),C^-_\Delta(\frk,\frl)]+[C(\frg,\frk),C^-(\frg,\frk)] = D_\Delta(\frk,\frl)^2+D(\frg,\frk)^2.
$$
This last equality can also be obtained from Kostant's formula for $D^2$:

\begin{multline*}
\Delta_{\frg,\frk}(D(\frk,\frl))^2+D(\frg,\frk)^2 = \\
\Delta_{\frg,\frk}(\Omega_\frk\otimes 1-\Delta_{\frk,\frl}(\Omega_{\frl})+||\rho_\frk||^2-||\rho_\frl||^2) + (\Omega_\frg\otimes 1 -\Delta_{\frg,\frk}(\Omega_{\frk}) +||\rho_\frg||^2-||\rho_\frk||^2) = \\ \Omega_\frg\otimes 1 -\Delta_{\frg,\frl}(\Omega_{\frl}) +||\rho_\frg||^2-||\rho_\frl||^2
= D(\frg,\frl)^2.
\end{multline*}

Here we have used the more precise notation for the various diagonal embeddings to avoid confusion; note the equality $\Delta_{\frg,\frk}(\Delta_{\frk,\frl}(\Omega_{\frl})) =\Delta_{\frg,\frl}(\Omega_{\frl})$.

In other words, we see that the $\frg\frl(1,1)$ corresponding to the pair $(\frg,\frl)$ sits
diagonally in the direct product of the two copies of $\frg\frl(1,1)$ described above.
Let us summarize the above discussion:

\begin{cor} Assume $(\frg,\frk)$ is a hermitian symmetric pair, let $\frqq=\frl\oplus\fru$ be
a $\theta$-stable parabolic subalgebra of $\frg$, and assume that $\frl\subset\frk$ and 
$\fru\supset\frp^+$. Then there are two supercommuting copies of $\frg\frl(1,1)$ inside
$U(\frg)\otimes C(\frs)$. One is spanned by $C_\Delta(\frk,\frl)$, $C_\Delta^-(\frk,\frl)$, 
$\Delta E(\frk,\frl)$ and $\Delta D(\frk,\frl)^2$ and the other is spanned by
$C(\frg,\frk)$, $C^-(\frg,\frk)$, $E(\frg,\frk)$, and $D(\frg,\frk)^2$. The diagonal of the
product of these two super subalgebras is the copy of $\frg\frl(1,1)$ spanned by
$C$, $C^-$, $E$ and $D^2$ from the end of Section 2.  \qed
\end{cor}

\section{Hodge decomposition for $\bar\fru$-cohomology in hermitian case}

In this section, $(\frg,\frk)$ is a hermitian symmetric pair, $\frqq=\frl\oplus\fru$ is
a $\theta$-stable parabolic subalgebra of $\frg$, and we assume $\frl\subset\frk$ and 
$\fru\supset\frp^+$.

Let $V$ be a unitary $(\frg,K)$-module,
and consider the form $\langle\;,\rangle$ on $V\otimes S$ introduced in Section 4. 
To apply the results of Section 4, we decompose 
$$
V\otimes S = V\otimes S_\frp\otimes S_{\frs\cap\frk}=
V\otimes {\textstyle\bigwedge^\cdot}\frp^+\otimes {\textstyle\bigwedge^\cdot}\fru\cap\frk,
$$
and embed $V\otimes \bigwedge^\cdot\frp^+$ as $V\otimes \bigwedge^\cdot\frp^+\otimes 1$. 
The form $\langle\;,\rangle$ restricts to the analogous definite form on 
$V\otimes \bigwedge^\cdot\frp^+$. 

Denote as before by $D=D(\frg,\frl)$ the Dirac operator for the pair $(\frg,\frl)$ and by
$C=C(\frg,\frl)$ and $C^-=C^-(\frg,\frl)$ its parts coresponding to $\fru$ and $\bar\fru$.
By Corollary 6.7, $C=C_\Delta(\frk,\frl)+C(\frg,\frk)$, and similarly for $C^-$. Moreover,
the copy of $\frg\frl(1,1)$ corresponding to the pair $(\frg,\frk)$ supercommutes with 
the copy of $\frg\frl(1,1)$ corresponding to the pair $(\frk,\frl)$.

By Corollary 4.3, the adjoints of $C_\Delta(\frk,\frl)$ and $C(\frg,\frk)$ are respectively
$C^-_\Delta(\frk,\frl)$ and $-C^-(\frg,\frk)$. So the adjoint of $C$ is 
$C^{adj}=C^-_\Delta(\frk,\frl) -C^-(\frg,\frk)$.

We consider the positive semidefinite operator $\Delta=CC^{adj}+C^{adj}C=[C,C^{adj}]$.
By the above remarks we have 

\begin{multline*}
\Delta=[C_\Delta(\frk,\frl)+C(\frg,\frk),C^-_\Delta(\frk,\frl)-C^-(\frg,\frk)]= \\
[C_\Delta(\frk,\frl),C^-_\Delta(\frk,\frl)]-[C(\frg,\frk),C^-(\frg,\frk)]=
D_\Delta(\frk,\frl)^2-D(\frg,\frk)^2.
\end{multline*}

We know from Section 4 that $V\otimes \bigwedge^\cdot\frp^+$ decomposes into eigenspaces of $D(\frg,\frk)^2$ for eigenvalues $\lambda\leq 0$.
Each eigenspace is $\tilde K$-invariant, and each $\tilde K$-isotypic component of 
$V\otimes \bigwedge^\cdot\frp^+$ is contained in an eigenspace. We assume $V$ is admissible,
so the eigenspaces are finite-dimensional.

Passing from $V\otimes \bigwedge^\cdot\frp^+$ to $V\otimes S$ is tensoring with the 
finite-dimensional $\frl$-module $\bigwedge^\cdot\fru\cap\frk$.
On this last space, there is no action of $U(\frg)$ or $U(\frk_\Delta)$. 
So every eigenspace of $D(\frg,\frk)^2$ on $V\otimes \bigwedge^\cdot\frp^+$ just gets tensored 
with $\bigwedge^\cdot\fru\cap\frk$, and this gives the eigenspace on $V\otimes S$ for the
same eigenvalue.

Since $D_\Delta(\frk,\frl)^2$ commutes with $D(\frg,\frk)^2$, it preserves these eigenspaces.
Moreover, the Levi subgroup $L\subset K$ corresponding to $\frl$ is compact. So is then the
double cover $\tilde L$, which acts on $V\otimes S$. Since $\tilde L$ commutes 
with $D(\frg,\frk)^2$, it also preserves its eigenspaces and hence these eigenspaces decompose
into $\tilde L$-irreducibles.
Since $D_\Delta(\frk,\frl)^2$ is up to a constant equal to the Casimir element of $\frl_\Delta$,
it follows that $D_\Delta(\frk,\frl)^2$ 
diagonalizes on each eigenspace of $D(\frg,\frk)^2$. To conclude:

\begin{lemma} $V\otimes S$ is a direct sum of eigenspaces for $\Delta$. In particular,
$V\otimes S=\ker\Delta\oplus\im\Delta$. \qed
\end{lemma}

This is an analogue of Corollary 4.5 for $\Delta$ in place of $D^2$. Now the arguments proving 
Lemma 4.6 and Theorem 4.7 work without change, and we obtain

\begin{thm} (a) $\ker\Delta = \ker C\cap\ker C^{adj}$;

(b) $V\otimes S = \ker\Delta\oplus\im C\oplus\im C^{adj}$;

(c) $\ker C=\ker\Delta\oplus\im C$;

(d) $\ker C^{adj}=\ker\Delta\oplus\im C^{adj}$.

\end{thm}

In other words, we have obtained a Hodge theorem for $\bar\fru$-cohomology.

To obtain it also for $\fru$-homology, we note that 
$(C^-)^{adj}=(C^-_\Delta(\frk,\frl)+C^-(\frg,\frk))^{adj}=
C_\Delta(\frk,\frl)-C(\frg,\frk)$, and so 

\begin{multline*}
[C^-,(C^-)^{adj}]=[C^-_\Delta(\frk,\frl)+C^-(\frg,\frk),C_\Delta(\frk,\frl)-C(\frg,\frk)]=\\
[C^-_\Delta(\frk,\frl),C_\Delta(\frk,\frl)]-[C^-(\frg,\frk),C(\frg,\frk)]=\Delta.
\end{multline*}

So the situation for $C^-$ is exactly the same as for $C$ and we conclude

\begin{thm} (a) $\ker\Delta = \ker C^-\cap\ker (C^-)^{adj}$;

(b) $V\otimes S = \ker\Delta\oplus\im C^-\oplus\im (C^-)^{adj}$;

(c) $\ker C^-=\ker\Delta\oplus\im C^-$;

(d) $\ker (C^-)^{adj}=\ker\Delta\oplus\im (C^-)^{adj}$.

\end{thm}

In other words, Hodge decomposition also holds for $\fru$-homology. Moreover, we see that
$\bar\fru$-cohomology and $\fru$-homology have the same set of harmonic representatives, 
$\ker\Delta$. In particular they are isomorphic.

We now want to relate $\bar\fru$-cohomology and $\fru$-homology  to Dirac cohomology with respect 
to $D(\frg,\frl)$. The main observation here is

\begin{lemma}
$\ker\Delta=\ker D_\Delta(\frk,\frl)^2\cap\ker D(\frg,\frk)^2=
\ker D_\Delta(\frk,\frl)\cap\ker D(\frg,\frk)$.
\end{lemma}

\begin{proof} The operators $D_\Delta(\frk,\frl)^2$ and $-D(\frg,\frk)^2$ are both positive
semidefinite and their sum is $\Delta$. This immediately implies the first equality.
The second follows from $\ker D_\Delta(\frk,\frl)^2=\ker D_\Delta(\frk,\frl)$ (since
$D_\Delta(\frk,\frl)$ is self-adjoint) and $\ker D(\frg,\frk)^2=\ker D(\frg,\frk)$
(since $D(\frg,\frk)$ is anti-self-adjoint). 
\end{proof}

We can now combine Theorems 7.2 and 7.3 with Lemmas 7.4 and 4.6 to conclude

\begin{cor} 
$\ker\Delta=\ker C_\Delta(\frk,\frl)\cap\ker C^-_\Delta(\frk,\frl) \cap\ker C(\frg,\frk)\cap\ker 
C^-(\frg,\frk).$ 
\end{cor}

\begin{proof} It is obvious that the left hand side contains the right hand side. Conversely, if
$x\in\ker\Delta$, then $Cx=0$ by Theorem 7.2, $D(\frg,\frk)x=0$ by Lemma 7.4, so 
$C(\frg,\frk)x=0$ by Lemma 4.6 and so also $Cx-C(\frg,\frk)x= C_\Delta(\frk,\frl)x=0$.
Analogously, $C^-(\frg,\frk)x=0$ and $C^-_\Delta(\frk,\frl)x=0$.
\end{proof}

Since $\ker C_\Delta(\frk,\frl)\cap\ker C(\frg,\frk)$ can be thought of as the kernel of
$C_\Delta(\frk,\frl)$ acting on the kernel of $C(\frg,\frk)$, and similarly for the 
$C^-$-operators, in view of Theorem 4.7 and Remark 4.8 we can reinterprete Corollary 7.5 
as follows: 

\begin{cor} To calculate the $\bar\fru$-cohomology of $V$, one can first calculate the
$\frp^-$-cohomology of $V$ to obtain a $\tilde K$-module, 
and then calculate the $\bar\fru\cap\frk$-cohomology of this module.
Analogously, to calculate the $\fru$-homology of $V$, 
one can first calculate the
$\frp^+$-homology of $V$, and then the $\fru\cap\frk$-homology of the resulting $\tilde K$-module.
\qed
\end{cor}

\begin{rmk}
Note that this is in fact the Hochschild-Serre spectral sequence
for the ideal $\frp^-$ of $\bar\fru$ respectively the ideal $\frp^+$ of $\fru$.
What we have obtained is that these Hochschild-Serre spectral sequences are always degenerate
for a unitary $(\frg,K)$-module $V$. \qed
\end{rmk}

We now turn our attention to the Dirac cohomology of $D=D(\frg,\frl)$. In addition to
the above considerations, we bring in
Corollary 5.6, and note that for both  $D_\Delta(\frk,\frl)$ and $D(\frg,\frk)$ the cohomology
is the same as the kernel or the kernel of the square. Thus we obtain:

\begin{thm} The Dirac cohomology $H_D(\frg,\frl;V)$ of a unitary $(\frg,K)$-module $V$ is 
isomorphic to the $\bar\fru$-cohomology of $V$ and the $\fru$-homology of $V$ up to
appropriate modular twists.
Moreover, all three cohomologies have the same set of harmonic representatives,
$\ker\Delta$. \qed
\end{thm}

\section{Homological properties of Dirac cohomology}

Let us start by showing that although we proved that in some cases Dirac cohomology of a unitary $(\frg,K)$-module with respect to $D(\frg,\frl)$ can
be identified with $\bar\fru$-cohomology or $\fru$-homology, one should by no means expect
that these notions agree for general $(\frg,K)$-modules. Let us see that this is not the case
even for $(\frs\frl(2,\bbC),SO(2))$-modules. 

Consider the module $V$ which is a nontrivial extension of the discrete series representation $W$
of highest weight $-2$ by the trivial module $\bbC$:
$$
0\to \bbC\to  V \to W\to 0.
$$
(In other words, $V$ is a dual Verma module.) 
The $\frk$-weights of $V$ (for the basis element 
$\bigl[\begin{smallmatrix} 0&-i\\ i&0\end{smallmatrix}\bigr]$ of $\frk$) are $0,-2,-4,\dots$. 
We are considering the case $\frl=\frk$, $\fru$ is spanned by 
$u=\frac{1}{2}\bigl[\begin{smallmatrix}  1&i\cr i&-1\end{smallmatrix}\bigr]$ and $\bar\fru$ is spanned by 
$u^*=\frac{1}{2}\bigl[\begin{smallmatrix} 1&-i\cr -i&-1\end{smallmatrix}\bigr]$.

For any $(\frs\frl(2,\bbC),SO(2))$-module $X$ we have
$$
X\otimes S = X\otimes 1\;\;\oplus\;\; X\otimes u,
$$
with $d:X\otimes 1\to X\otimes u$ given by $d(v\otimes 1)=u^*\cdot v\otimes u$, 
$\del:X\otimes u\to X\otimes 1$ given by $\del(v\otimes u)=u\cdot v\otimes 1$, and $D=d-2\del$.
By an easy direct calculation, we see that

\begin{gather*}
H_0(\fru;V)=0;\qquad H_1(\fru;V)=\bbC v_0\otimes u; \\
H^0(\bar\fru;V)=\bbC v_0\otimes 1;\qquad H^1(\bar\fru;V)=\bbC v_0\otimes u\oplus \bbC v_{-2}\otimes u; \\
H_D(V)=\Ker D = \bbC v_0\otimes u.
\end{gather*}
as vector spaces.
Here $v_i$ denotes a vector in $V$ of $\frk$-weight $i$. So we see
$$
H_D(V)=H_\cdot(\fru;V)\neq H^\cdot(\bar\fru;V)
$$

On the other hand, for $\bbC$ and $W$ another easy calculation
(which can be shortened by using Theorem 4.7, as both $\bbC$ and $W$ are unitary)
implies

\begin{gather*} 
H^0(\bar\fru;\bbC)=H_0(\fru;\bbC)=H_D^0(\bbC)=\bbC 1\otimes 1; \\
H^1(\bar\fru;\bbC)=H_1(\fru;\bbC)=H_D^1(\bbC)=\bbC 1\otimes u;\\
H^0(\bar\fru;W)=H_0(\fru;W)=H_D^0(W)=0; \\
H^1(\bar\fru;W)=H_1(\fru;W)=H_D^1(W)=\bbC v_{-2}\otimes u.
\end{gather*}

To explain why Dirac cohomology of $V$ differs from the $\bar\fru$-cohomology
of $V$, we will examine their behavior with respect to extensions.

Recall the well known long exact sequences for Lie algebra homology and
cohomology corresponding to our short exact sequence
$$
0\to \bbC\to V\to W\to 0.
$$
They are
$$
0\to \bbC 1\otimes u\to \bbC v_0\otimes u \overset{0}\to  \bbC v_{-2}\otimes u \to \bbC 1\otimes 1 \to 0\to 0\to 0
$$
for $\fru$-homology and
$$
0\to  \bbC 1\otimes 1 \to  \bbC v_0\otimes 1\to  0 \to \bbC 1\otimes u
\to \bbC v_{-2}\otimes u\oplus\bbC v_0\otimes u \to  \bbC v_{-2}\otimes u \to  0
$$
for $\bar\fru$-cohomology. Here all arrows are the obvious ones except for the
one labelled by $0$.

For Dirac cohomology, instead of a long exact sequence (which clearly does
not makes sense in general, as Dirac cohomology is not $\bbZ$-graded), there
is a six-term exact sequence. In the above example, this sequence is
$$
\CD
\bbC 1\otimes 1 @>>> 0 @>>> 0 \\
@AAA @.                   @VVV  \\
\bbC v_{-2}\otimes u @<0<<\bbC v_0\otimes u @<<< \bbC 1\otimes u
\endCD
$$
We see that in this example the six-term sequence agrees with the $\fru$-homology long exact
sequence only because of the presence of zeros; in general, all three sequences are different.

To define the six-term sequence in a more general situation, let us assume that the
Dirac cohomology is $\bbZ_2$-graded. This happens
whenever we are starting from $\frg=\frr\oplus\frs$
with $\frs$ even-dimensional; this is automatic when $\frr=\frl$ is a Levi subalgebra. 
Let 
$$
0\to X \overset{i}\to Y\overset{p}\to Z\to 0
$$
be a short exact sequence of $(\frg,K)$-modules. Tensor this sequence by $S$, and denote
the arrows still by $i$ and $p$ (they get tensored by the identity on $S$).
Assuming that $D^2$ is a semisimple operator
for each of the three modules, we can construct a six-term exact sequence
$$
\CD
H_D^0(X)@>>>H_D^0(Y)@>>>H_D^0(Z) \\
@AAA        @.          @VVV     \\
H_D^1(Z)@<<<H_D^1(Y)@<<<H_D^1(X) 
\endCD
$$
The horizontal arrows are induced by $i$ and $p$. The vertical arrows are the
connecting homomorphisms, defined as follows. Let $z\in Z\otimes S$ represent a Dirac cohomology
class, so $Dz=0$. Choose $y\in Y\otimes S$ such that $p y=z$. Since $D^2$ is semisimple, we
can assume $D^2y=0$. Since $pDy=Dpy=Dz=0$, we see that $Dy=ix$ for some $x\in X$. Since
$D^2y=0$, we see that $Dx=0$, so $x$ defines a cohomology class. This class is by definition
the image of the class of $z$ under the connecting homomorphism. Clearly, we changed parity when
we applied $D$, and this defines both vertical arrows at once.

It is easy to see that this map is well defined, and that the resulting six-term sequence is exact.
To conclude:

\begin{thm} Let $\frg=\frr\oplus\frs$ be an orthogonal decomposition, with $\frr$ a
reductive subalgebra and $\frs$ even-dimensional. Let $0\to X \to Y\to Z\to 0$ be a 
short exact sequence of $(\frg,K)$-modules and
assume that the square of the Dirac operator $D(\frg,\frr)$
is a semisimple operator for $X$, $Y$ and $Z$. Then there is a six-term exact sequence 
corresponding to this short exact sequence, as described above. \qed
\end{thm} 

Finally, let us comment on what can be done when $\frs$ is odd dimensional, say $\dim\frs=2n+1$.
The usual spin modules $S_1$ and $S_2$ are not $\bbZ_2$-graded and thus it seems the above construction does not make sense. Recall that $S_1$ and $S_2$ are defined by
writing $\frs=\bbC^{2n}\oplus\bbC$, considering the spin module for
$C(\bbC^{2n})$, and letting the last basis element of $\frs$ act in two
different ways (preserving the even and odd subspace instead of exchanging them).

We can instead consider the unique irreducible graded module $\tilde S$ of $C(\frs)$. It
can be constructed as the restriction of the (unique) spin module for $C(\bbC^{2n+2})$ to
$C(\frs)\subset C(\bbC^{2n+2})$. As a non-graded module, $\tilde S$ decomposes as $S_1\oplus S_2$.
If we define Dirac cohomology using $\tilde S$ in place of $S_1$ or $S_2$, we double it, but we
do get a $\bbZ_2$-grading. Then the above construction works also in the odd case. Thus, this
is probably a more natural definition of Dirac cohomology in the odd case.

It remains to see what can be done if $D^2$ is not a semisimple operator. One possibility
might be to consider a more general definition of Dirac cohomology in that setting.


\begin{thebibliography}{999}

\bibitem [A]{A} J.~Adams, 
\emph{Nilpotent cohomology of the oscillator representation},
 J. Reine Angew. Math. \textbf{449} (1994), 1--7. 

\bibitem [AM]{AM} A.~Alekseev, E.~Meinrenken,
\emph{Lie theory and the Chern-Weil homomorphism}, math.RT/0308135.

\bibitem [AS]{AS} M.~Atiyah, W.~Schmid,
\emph{A geometric construction of the discrete series
for semisimple Lie groups},
Invent. Math. \textbf{42} (1977), 1--62. 

\bibitem [BW]{BW} A.~Borel, N.~Wallach,  
\emph{Continuous cohomology, discrete subgroups, and representations of reductive groups}, 
Second edition. Mathematical Surveys and Monographs, 67. American Mathematical Society, Providence, RI, 2000.

\bibitem [CO]{CO}  W.~Casselman, M.~S.~Osborne,
\emph{The $\frn$-cohomology of representations with an  infinitesimal character}, 
Comp. Math. \textbf{31} (1975), 219--227.

\bibitem [Ch]{Ch} C.~Chevalley,
\emph{The algebraic theory of spinors},
Columbia University Press, 1954.

\bibitem [C]{C} D.~H.~Collingwood,
\emph{The $\frn$-homology of Harish-Chandra modules: generalizing a theorem of Kostant}, 
Math. Ann. \textbf{272} (1985), no. 2, 161--187. 

\bibitem [CM]{CM} A.~Connes, H.~Moscovici,
\emph{The $L^2$-index theorem for homogeneous spaces of Lie groups},
Ann. of Math. \textbf{115} (1982), no. 2, 291--330.

\bibitem [E]{E} T.~J.~Enright, 
\emph{Analogues of Kostant's $\fru$-cohomology formulas for unitary highest weight modules}, J. Reine Angew. Math. \textbf{392} (1988), 27--36.

\bibitem [HoP]{HoP} R.~Hotta, R.~Parthasarathy,
\emph{A geometric meaning of the multiplicities of integrable discrete classes in $L^2(\Gamma\backslash G)$},
Osaka J. Math. \textbf{10} (1973), 211--234.

\bibitem [HP1]{HP1} J.-S.~Huang, P.~Pand\v zi\'c,
\emph{Dirac cohomology, unitary
representations and a proof of a conjecture of Vogan},
J. Amer. Math. Soc. \textbf{15} (2002), 185--202.

\bibitem [HP2]{HP2} J.-S.~Huang, P.~Pand\v zi\'c,
\emph{Dirac operators in representation theory}, 
Representations of Real and P-adic Groups,
Lecture Notes Series, Vol. 2, 
Institute for Mathematical Sciences, National University of Singapore,
Singapore University Press and World Scientific, 2004,
pp. 163-219.

\bibitem [Kac]{Kac} V.~Kac, \emph{Lie superalgebras},
Adv. in Math. \textbf{26} (1977), 8--96.

\bibitem [K1]{K1} B.~Kostant,
\emph{Lie algebra cohomology and the generalized Borel-Weil theorem},
Ann. of Math. \textbf{74} (1961), 329--387.

\bibitem [K2]{Kcubic} B.~Kostant, 
\emph{A cubic Dirac operator and the emergence of
Euler number multiplets of representations for equal rank subgroups},
Duke Math. Jour. \textbf{100}  (1999), 447--501.

\bibitem [K3]{Kgen} B.~Kostant,
\emph{A generalization of the Bott-Borel-Weil theorem and Euler number multiplets of representations}, Lett. Math. Phys. \textbf{52} (2000), 61--78. 

\bibitem [K4]{Kdircoh} B.~Kostant,
\emph{Dirac cohomology for the cubic Dirac operator},
Studies in memory of
I. Schur, in Progress of Math. vol. \textbf{210} (2003), 69--93.

\bibitem [Ku]{Ku} S.~Kumar,
\emph{Induction functor in non-commutative equivariant cohomology and Dirac cohomology},
preprint, University of North Carolina, 2003.

\bibitem [L]{L} R.~Langlands,
\emph{The dimension of spaces of automorphic forms},
Amer. J. of Math. \textbf{85} (1963), 99--125.

\bibitem [P]{P}  R.~Parthasarathy
\emph{Dirac operator and the discrete series},
Ann. of Math. \textbf{96} (1972), 1--30. 

\bibitem [S]{S} W.~Schmid, 
\emph{$L^2$-cohomology and the discrete series},
Ann. of Math. \textbf{103} (1976), 375--394.

\bibitem [V1]{V1} D.A.~Vogan,
\emph{Representations of real reductive Lie groups},
Birkh\"auser, Boston, 1981.

\bibitem[V2]{V2}
D.A.~Vogan, \emph{Dirac operators and unitary representations},
3 talks at MIT Lie groups seminar, Fall 1997.

\bibitem[V3]{V3}
D.A.~Vogan, \emph{$\mathfrak n$-cohomology in representation theory},
a talk at ``Functional Analysis VII'', Dubrovnik, Croatia, September 2001.

\bibitem[W]{W} N.~R.~Wallach,
\emph{Real Reductive Groups, Volume I},
Academic Press, 1988.


\end{thebibliography}
\end{document}